\documentclass[12pt]{article}
\usepackage{amsmath,amssymb,amsthm,setspace}
       \newtheorem{theorem}{Theorem}
       
       \newtheorem{lemma}[theorem]{Lemma}

       {\theoremstyle{definition}

       \newtheorem{conjecture}{Conjecture}}
\def\N{\ensuremath{\mathbb N}} 
\def\R{\ensuremath{\mathbb R}} 
\def\C{\ensuremath{\mathbb C}} 
\def\Z{\ensuremath{\mathbb Z}} 

\newcommand{\can}{{\cal N}}

\newcommand{\lf}{\lfloor}
\newcommand{\rf}{\rfloor}
\newcommand{\tlll}{{\widetilde{\lambda}}}
\newcommand{\tX}{{\widetilde{X}}}
\newcommand{\tA}{{\widetilde{A}}}
\newcommand{\tSSS}{{\widetilde{\Sigma}}}
\newcommand{\tmmm}{{\widetilde{\mu}}}
\newcommand{\tT}{{\widetilde{T}}}
\newcommand{\tx}{{\widetilde{x}}}
\newcommand{\tf}{{\widetilde{f}}}
\newcommand{\off}{{\overline{f}}}
\newcommand{\tq}{{\widetilde{q}}}

\newcommand{\cae}{{\cal E}}

\newcommand{\defeq}{{\buildrel {\rm def}\over =}}

\newcommand{\tqjm}{{\widetilde{q}}_{j,m}}

\newcommand{\oN}{\overline{N}}
\newcommand{\offf}{\overline{\phi}}
\newcommand{\fff}{{\phi}}

\newcommand{\hfff}{\widehat{\phi}}
\newcommand{\vppp}{\check{\psi}}
\newcommand{\cab}{{\cal B}}

\newcommand{\eee}{\epsilon}

\newcommand{\aaa}{\alpha}
\newcommand{\bbb}{\beta}

\newcommand{\oS}{{\overline{S}}}

\newcommand{\LLL}{\Lambda}
\renewcommand{\lll}{\lambda}

\renewcommand{\ggg}{\gamma}

\newcommand{\kkk}{\kappa}
\newcommand{\oo}{\infty}
\newcommand{\sse}{\subset}
\newcommand{\ess}{\emptyset}
\newcommand{\sm}{\setminus}

\title{Universally $L^1$ good sequences with gaps tending to
infinity}
\author{
Zolt\'an Buczolich\thanks{
Research supported by the Hungarian
National Foundation for Scientific research T049727.
},
Department of Analysis, E\"otv\"os Lor\'and\\
University, P\'azm\'any P\'eter S\'et\'any 1/c, 1117 Budapest, Hungary\\
email: buczo@cs.elte.hu\\
{\tt www.cs.elte.hu/\hbox{$\sim$}buczo}
}
\date{\today}
\begin{document}
\maketitle
\begin{abstract}
We construct a sequence $(n_{k})$ such that
$n_{k+1}-n_{k}\to\oo$ and for any ergodic dynamical
system $(X,\Sigma,\mu,T)$ and $f\in L^{1}(\mu)$
the averages $\lim_{N\to\oo}(1/N)\sum_{k=1}^{N}f(T^{n_{k}}x)$
converge to $\int_{X}fd\mu$ for $\mu$ almost every $x.$
Since the above sequence is of zero Banach density
this disproves a conjecture of
J. Rosenblatt and M. Wierdl about the nonexistence
of such sequences.
\end{abstract}

\section{Introduction}

In \cite{[DSA]} it is shown that the sequence $n_{k}=k^{2}$,
$k=1,2,...$ is $L^{1}$-universally bad. This means that for
all aperiodic ergodic dynamical systems $(X,\Sigma,\mu,T)$ there
exists $f\in L^{1}(\mu)$ such that the ergodic averages
\begin{equation}\label{*I1}
\lim_{N\to\oo}\sum_{k=1}^{N}f(T^{k^{2}}x)
\end{equation}
fail to converge on a set of positive measure.
On the other hand, results of Bourgain
\cite{B1}, \cite{B2} and \cite{B3},
imply that \eqref{*I1} converges $\mu$ almost everywhere
whenever $f\in L^{p}(\mu)$ with $p>1$.

When I was working on \cite{[DSA]} I learned from M. Keane
that it is not known whether there exists a sequence $(n_{k})$
such that $n_{k+1}-n_{k}\to\oo$ and for any $f\in L^{1}(\mu)$
\begin{equation}\label{*I2}
\lim_{N\to\oo}\frac{1}N\sum_{k=1}^{N} f(T^{n_{k}}x)
\end{equation}
converges $\mu$ almost everywhere.
This question is also stated in \cite{[Jo]}
on p. 64 in the second paragraph of Section 7.
A sequence satisfying
$n_{k+1}-n_{k}\to\oo$ is of zero Banach density.
In \cite{[RW]} the authors ``risk" the following conjecture
(see Conjecture 4.1 on p. 74 of \cite{[RW]}, here we use
slightly different equivalent notation):
\begin{conjecture}\label{conjrw}
Suppose that the sequence $(n_{k})$ has zero Banach density and
let $(X,\Sigma,\mu, T)$ be an aperiodic dynamical system.
Then for some $f\in L^{1}(\mu)$ the averages
\eqref{*I2} do not converge almost everywhere.
\end{conjecture}

The purpose of this paper is to show that that there exist
universally $L^{1}$-good sequences $(n_{k})$ for which
$n_{k+1}-n_{k}\to\oo.$ A sequence is universally
$L^{1}$-good if \eqref{*I2} converges $\mu$ almost everywhere
for any ergodic dynamical system $(X,\Sigma,\mu,T)$
and $f\in L^{1}(\mu).$ This implies that Conjecture \ref{conjrw}
is false. This also provides an explanation why was it
so difficult to obtain the result in \cite{[DSA]} which states
that $n_{k}=k^2$ is $L^{1}$-universally bad.

In this paper, like in \cite{B1}, we mean by a dynamical
system $(X,\Sigma,\mu,T)$ an invertible measure preserving
 transformation acting on a probability measure space. We also
assume that $T$ is aperiodic.
By scrutinizing the proof presented in this paper one can see
that for our sequence $(n_{k})$
the averages \eqref{*I2} converge almost everywhere in
ergodic periodic systems as well. The non-invertible
case from the point of view of this paper can easily
be
reduced to the invertible one. Suppose that for a
non-invertible aperiodic ergodic dynamical system
$(X,\Sigma,\mu,T)$ there exists $f\in L^{1}(\mu)$
for which \eqref{*I2} diverges when $x\in A\in \Sigma$
and $\mu(A)>0.$ Consider the natural extension
$(\tX,\tSSS,\tmmm,\tT)$ of $(X,\Sigma,\mu,T)$
(see \cite{[CFS]}, Chapter 10, \S 4., or \cite{[P]}
1.3.G., p. 13).
Then $(\tX,\tSSS,\tmmm)$ is the inverse limit space
obtained from $(X,\Sigma,\mu,T)$.
The elements of $\tX$ are of the form
$\tx=(x_{0},x_{1},...)$ with $T(x_{j})=x_{j-1}$,
$j=1,2,....$
The transformation $\tT$ is defined so that
$\tT\tx=(Tx_{0},Tx_{1},...).$
Then $\tT^{-1}\tx=(x_1,x_{2},...)$ and by
 Theorem 1, on p. 241 of \cite{[CFS]}
$\tT$ is an ergodic measure preserving transformation.

 Set $\tA= \{ \tx\in\tX:
x_{0}\in A  \}$. Then $\tmmm(\tA)=\mu(A)>0$.
If we set $\tf(\tx)=f(x_{0})$ then $\tf\in
L^{1}(\tmmm)$ and $$\lim_{N\to\oo}(1/N)\sum_{k=1}^{N}
\tf(\tT^{n_{k}}\tx)
=\lim_{N\to\oo}(1/N)\sum_{k=1}^{N}
f(T^{n_{k}}x_{0})$$ diverges for all $\tx\in\tA.$
This shows that if $(n_{k})$ is $L^{1}$-bad for
a non-invertible system then it is also
bad for a suitable invertible one.

This paper is organized as follows. After this introduction
in Section \ref{main} we state Theorem \ref{mainth}
which is the main result of this paper about the existence
of universally $L^{1}$-good sequences $(n_{k})$
with gaps converging to infinity. The proof of
Theorem \ref{mainth} is based on Lemmas
\ref{dl1} and \ref{*dl2}. In Lemma \ref{dl1}
we show that the $(n_{k})$ averages converge
for simple functions, which form a dense
subset in $L^{1}$. In Lemma \ref{*dl2} a weak
$(1,1)$ inequality is established for the maximal
operator corresponding to the sequence $(n_{k})$.

In Section \ref{defnk} we define $(n_{k})$ by induction.
Intervals $[\bbb_{m-1},\bbb_{m})$ are selected and
at the $m$'th step of our definition we choose the
terms of $(n_{k})$ in one such interval. One can think
of the terms of $(n_{k})$ as the union of
finitely many arithmetic sequences with those
terms deleted which are too close to each other.
To be more specific, we choose $K_{m}$ many
different prime numbers $q_{j,m}$ and consider
those terms of the set $ \{ lq_{j,m}:l\in\Z, j=1,...,K_{m}  \}$
which are
in $[\bbb_{m-1},\bbb_{m})$ and delete those ones which
are too close.

In Section \ref{secfunz} we consider functions on
$\Z$ with bounded support. We introduce the operators $\cab$
and $\cab_{0}$ with maximal operators $\cab^{*}$ and
$\cab_{0}^{*}$. The maximal inequalities established
in this section will be applied in later sections
 with a fixed $m\in \Z$ for the terms of $(n_{k})$ which are
in $[\bbb_{m-1},\bbb_{m}).$ The most important result is
in Lemma \ref{d8} about $\cab_{0}^{*}$.
Lemmas \ref{Lemma3.5} and \ref{*strong}
are mere restatements of well-known maximal inequalities.
Lemma \ref{dl24} contains a not too difficult maximal
inequality about the operator $\cab^{*}$.

In Sections \ref{secpart1} and \ref{secpart2} we prove Lemma
\ref{*dl2}. The second part of the proof of Lemma 3, given
in Section \ref{secpart2} is used for the proof of
Lemma \ref{dl1} as well.
This means that some estimates and notation introduced here
is used only later, in Section \ref{sec7}. This shared proof
part explains that instead of using some kind of transference
principle why we use directly
 Kakutani-Rokhlin tower constructions in Sections
\ref{secpart1} and \ref{secpart2} to transfer the
results from Section \ref{secfunz} to arbitrary
dynamical systems. Of course, we also need to
``paste" together the estimates which we obtain
for different $m$'s for terms of $(n_{k})$
in $[\bbb_{m-1},\bbb_{m}).$
To estimate the $(n_{k})$ averages of \eqref{*I2}
we represent $f$ as $f=\lll'(f_{1,m}+f_{2,m}+f_{3,m})$ with
$\lll'\in\R$ and $m\in\N.$
In Section \ref{secpart1} we deal with terms
involving $f_{2,m}$ and $f_{3,m}$. While
the terms involving $f_{1,m}$ are estimated in Section
\ref{secpart2}.
If $f$ is bounded and $N$ is large then in \eqref{*I2}
we can replace $f$ by $\lll'f_{1,m}$ and this is
why Section \ref{secpart2} is used in the proof
of Lemma \ref{dl1} as well.
In Section \ref{secpart1} during the estimates
related to the terms $f_{2,m}$ an operator
denoted by $B$ is defined. After the
Kakutani-Rokhlin tower construction it
turns out that $B$ coincides with $\cab$
and the simpler maximal inequality of
Lemma \ref{dl24} can be used to estimate
the maximal operators $B^{*}$ and
$\cab^{*}$. It simplifies our work that by
\eqref{*d18b}, $\sum_{m}f_{2,m}\leq 3f/\lll'$
and hence $\sum_{m}f_{2,m}\in L^{1}.$
Unfortunately, it is not always true that
$\sum_{m} f_{1,m}\in L^{1}.$ This is why
we need in Section \ref{secpart2} much
more sophisticated methods than the ones in
Section \ref{secpart1}. Here we need to introduce
the modified operators $B_{0}$ which
coincide with $\cab_{0}$ after the
Kakutani-Rokhlin tower construction.
In this section the more involved
Lemma \ref{d8} is needed
for the estimation of the maximal
operators  $B_{0}^{*}$ and
$\cab_{0}^{*}$.

In Section \ref{sec7} based on Part 2
of the proof of Lemma \ref{*dl2} we see
that for simple functions the
$(n_{k})$-averages in \eqref{*I2}
do not differ much from the ordinary
ergodic averages and hence Birkhoff's
Ergodic theorem implies Lemma \ref{dl1}.

\section{Main Result}\label{main}

The desired universally $L^{1}$-good sequence
with gaps tending to infinity
will be denoted by $( n_{k})$.

We set
$$\oN_{a}^{b}=\# \{ n_{k}:n_{k}\in [a,b)  \}.$$

Suppose
$f\in L^{1}(\mu)$.
We set
$$A(f,x,N)=\frac{1}{\oN_{0}^{N}}
\sum_{k=1}^{\oN_{0}^{N}}f(T^{n_{k}}x).$$
We also introduce
$$A^{*}(f,x)=\sup_{1\leq N}|A(f,x,N)|.$$
The main result of the paper is the following:
\begin{theorem}\label{mainth}
There exists a sequence
$(n_{k})$ satisfying $n_{k+1}-n_{k}\to\oo$
(and hence of zero Banach density)
which is universally $L^{1}$-good, that is,
for any
invertible aperiodic ergodic dynamical system $(X,\Sigma,\mu,T)$
and $f\in L^{1}(\mu)$
we have
\begin{equation}\label{*main}
\lim_{N\to\oo}A(f,x,N)=\lim_{K\to\oo}\frac{1}K
\sum_{k=1}^{K}f(T^{n_{k}}x)= \int_{X}f d\mu,
\end{equation}
for $\mu$
almost every $x\in X.$
\end{theorem}

The proof of Theorem \ref{mainth} follows from the
following two lemmas.
The first one yields a dense set in $L^{1}$
for which the $A(f,x,N)$ averages converge.
A function $f:X\to\R$ is a simple function if it is
measurable and its range consists of a finite set.
\begin{lemma}\label{dl1}
With the assumptions of Theorem \ref{mainth},
for any simple function $f$ we have
\begin{equation}\label{**d2}
\lim_{N\to\oo}A(f,x,N)=\int_{X}f d\mu.
\end{equation}
\end{lemma}
The second one gives a weak $(1,1)$ inequality
for the maximal operator $A^{*}$.
\begin{lemma}\label{*dl2}
With the notation used in Theorem \ref{mainth}
for any $f\in L^{1}(\mu)$  and $\lll>0$
we have
\begin{equation}\label{*d2}
\mu( \{ x: A^{*}(f,x)>\lll  \})\leq \frac{1000||f||_{1}}
{\lll}.
\end{equation}
\end{lemma}

\begin{proof}[Proof of Theorem \ref{*main}]
By Lemma \ref{dl1} there exists a dense set of functions
in $L^{1}(\mu)$ for which
$\lim_{N\to\oo}A(f,x,N)=\int_{X}fd\mu$ holds $\mu$ almost everywhere.
The weak $(1,1)$ inequality of Lemma \ref{*dl2} then implies
the almost everywhere finiteness of the maximal
operator $A^{*}(f,x).$ By Banach's principle the
almost everywhere convergence of $A(f,x,N)$
follows for all $f\in L^{1}(\mu)$ (for the details
see \cite{[P]} 3.2., p. 91).
\end{proof}

For ease of notation, if we write $\int fd\mu$ we always mean
$\int_{X}fd\mu.$

\section{Definition of $(n_{k})$ and some estimates}\label{defnk}

We will use some intervals determined by the integers
$\bbb_{m}$. We set $\bbb_{-1}=\bbb_{0}=0$ and
the positive integers $\bbb_{1}<...<\bbb_{m}<...$ will be defined
by induction.
In each block
we will use different numbers $q_{j,m}$, $j=1,...,K_{m}$.
These numbers will be different primes if $m>1.$
Their product $p_{m}=q_{1,m}\cdot \cdot \cdot q_{K_{m},m}$
will be called the period used in block $m$.
We suppose that the primes $q_{j,m}$ are approximately
the same size, that is,
\begin{equation}\label{*d20f1}
\frac{1}2
< \frac{q_{j,m}}{q_{j',m}}<2\text{ for }j,j'\in \{ 1,...,K_{m}  \}.
\end{equation}
We put
$$\tqjm\defeq\frac{p_{m}}{q_{j,m}}, \text{ and }
Q(m)\defeq\sum_{j=1}^{K_{m}}\frac{1}{q_{j,m}}.$$
We will use a parameter $d_{m}$ which will be
a lower bound on the
distance among the
terms
of $(n_{k})$
belonging to the interval $[\bbb_{m-1},\bbb_{m}).$
We suppose that $d_{m}\to\oo$ and $d_{m}<q_{j,m}$
for all $j=1,...,K_m$. For example, the choice $d_{m}=m$
is suitable. The sequence $(d_{m})$ will ensure that
the gaps
between consecutive terms of
$(n_{k})$ converge to infinity and hence
$(n_{k})$ will have zero Banach density.

We put
$\oN_{-2}=\oN_{-1}=\oN_{0}=0$
and $$\oN_{m}=\oN_{0}^{\bbb_{m}}=
\# \{n_{k}: n_{k}\in [0,\bbb_{m})  \}.$$

We will choose our parameters so that $\oN_{m-1}$
is much larger than $p_{m}$ for $m=2,....$

Next we give the general plan of the definition
of our parameters by mathematical induction.
There will be several technical assumptions
about these parameters introduced later.
Here we just want to orientate the reader about
what is chosen, when.
To start our induction we put
$ K_{1}=1, q_{1,1}=1.$
Then
$p_{1}=1$ and
$Q(1)=1.$
At the first step, after $\bbb_{1}>10$
is determined, we will choose the terms of
$(n_{k})$ in $[\bbb_{0},\bbb_{1})$
so that $n_{k}=k-1$, for $k=1,...,\bbb_{1},$
that is, each integer from $[\bbb_{0},\bbb_{1})$
will belong to $(n_{k}).$

Suppose
for an $m>1$ we have
$\bbb_{m'-1}$, $K_{m'}$, and $q_{m',j}$
$j=1,...,K_{m'}$
for $m'\leq m-1$
and the terms of the sequence
$n_{k}$ which satisfy $n_{k}<\bbb_{m-2}$
are defined. This gives the values of
$\oN_{m'}$ for $m'\leq m-2$
as well.
Choose $K_{m}$ so that
\begin{equation}\label{*f13}
\frac{32}{K_{m}}\oN_{m-2}10^{4}\cdot 4^{m+1}<
2^{-(m+1)}.
\end{equation}
Next, one needs to choose  the prime numbers $q_{j,m}>d_{m}$
so that \eqref{*d20f1} holds, $p_{m}=q_{1,m}\cdot \cdot \cdot q_{K_{m},m}
>p_{m-1}$, $Q(m)<Q(m-1)$
and
\begin{equation}\label{*f15}
\oN_{m-2}\cdot 4\cdot K_{m}^{2}\cdot\frac{ d_{m}+1}
{\min_{j'} \{ q_{j',m}  \}}
<
\frac{1}{200 (m+1)}.
\end{equation}
For $m>3$ we also set
\begin{equation}\label{*f18}
\ggg_{m}=\frac{1}{2000\cdot (m+1)\cdot \oN_{m-2}}
<\frac{1}{2000\cdot m\cdot \oN_{m-3}},
\text{ and }\ggg_{\bbb }=\frac{1}{1000}.
\end{equation}
We put $\ggg_{1}
=\ggg_{2}=
\ggg_{3}=
\frac{1}8.$

After the selection of $p_{m}$ we choose
a sufficiently large
$\bbb_{m-1}$.

Later we need for $m=2,3,...$ that by our assumptions
\begin{equation}\label{*f4c}
1-\ggg_{m}>\frac{3}4\text{ and }p_{m-1}<p_{m}<\frac{1}{10^{4}}(\bbb_{m-1}-
\bbb_{m-2})
<\frac{1}{10^{4}}\bbb_{m-1}.
\end{equation}

The value of $\bbb_{m-1}$,
and the numbers $q_{j,m-1}$, $j=1,...,K_{m-1}$
 will determine the terms of $(n_{k})$
in $[\bbb_{m-2},\bbb_{m-1}).$ This will give
us the value of $\oN_{m-1}$ as well.
We will have several assumptions later about
$\bbb_{m-1}$ and $\oN_{m-1}$. One should think of
these assumptions that they require that
these numbers are
 much larger than similar parameters with
lower indices. For example, we will need that
\begin{equation}\label{*f19}
(\sum_{m'=1}^{m-2}\oN_{m'})\frac{\oN_{m-3}}{\oN_{m-1}}
<\frac{1}{3m}\text{ and }
p_{m}<\frac{1}{100}\oN_{m-1}<\frac{1}{100}\bbb_{m-1}.
\end{equation}
In addition,
for convenience,
we also suppose that
\begin{equation}\label{*pmbbb}
p_{m} \text{ divides }
\bbb_{m-1}.
\end{equation}

For ease of notation suppose that $\bbb_{m}$ and
the numbers $q_{j,m}$, $j=1,...,K_{m}$
are given for an $m>2.$
Next we discuss how these numbers determine $(n_{k})$
in $[\bbb_{m-1},\bbb_{m})$ for $m>1$. According to
\eqref{*pmbbb},
$p_{m}$ and hence all $q_{j,m}$ divide $\bbb_{m-1}$.
Set
$$\LLL_{j,m,0}=
 \{ lq_{j,m}:l\in \Z  \}\cap [\bbb_{m-1},\bbb_{m}).$$
If we take a union of the sets $\LLL_{j,m,0}$ for
$j=1,...,K_{m}$
then
some elements might be closer than $d_{m}$. So we need
to remove these points. First set
$$
\LLL_{j,m,0}^{-}= \bigg \{ n\in \LLL_{j,m,0}:
\exists\, n'\in\bigcup_{j'=1,\ j'\not=j}^{K_{m}}
\LLL_{j',m,0},\  |n'-n|\leq d_{m} \bigg  \}
$$
then put
$$\LLL_{m}=\bigcup_{j=1}^{K_{m}}\LLL_{j,m,0}\sm \LLL_{j,m,0}^{-}.$$
Since $\bbb_{m-1}$ belongs to all $\LLL_{j,m,0}$
we have
\begin{equation}\label{*spacing}
[\bbb_{m-1},\bbb_{m-1}+d_{m})\cap \LLL_{m}=\ess.
\end{equation}
We define the terms of $(n_{k})$ so that
$n_{k-1}<n_{k}$ and  $\{ n_{k}  \}_{k=1}^{\oo}\cap
[\bbb_{m-1},\bbb_{m})=\LLL_{m}
\cap [\bbb_{m-1},\bbb_{m}).$
Therefore,
letting $\LLL_{1}=[\bbb_{0},\bbb_{1})$
we have
$ \{ n_{k}:k=1,...  \}=\cup_{m\in \N}\LLL_{m}$
and the spacing of at least $d_{m}$ among the elements
of each $\LLL_{m}$ plus \eqref{*spacing} ensures that
$n_{k+1}-n_{k}\to\oo$ as $k\to\oo$.

Suppose $[n',n'+p_{m})\sse [\bbb_{m-1},\bbb_{m}).$
Then $$\#(\LLL_{j,m,0}\cap[n',n'+p_{m}))=\frac{p_{m}}
{q_{j,m}}=\tqjm,$$
and
\begin{equation}\label{*d4a}
\#(\LLL_{m}\cap [n',n'+p_{m}))\leq
\sum_{j=1}^{K_{m}}\#(\LLL_{j,m,0}\cap[n',n'+p_{m}))=
p_{m}Q(m).
\end{equation}
For $j'\not=j$ set
$$
\LLL_{j,m,0,j'}^{-}= \{
lq_{j,m}\in \LLL_{j,m,0}:\exists l'\in\Z,
\text{ such that }|lq_{j,m}-l'q_{j',m}|\leq d_{m}   \}.
$$

If $j'\not=j$
then $q_{j,m}$ and $q_{j',m}$ are relatively prime.
Modulo $q_{j',m}$ the numbers $lq_{j,m}$, $l=0,...,q_{j',m}-1$
hit each residue class exactly once. Hence, out of these
$2d_{m}+1$ are not farther than $d_{m}$
from $0$ modulo $q_{j',m}$.
Thus,
for each $j'\not=j$
 out of the $\tqjm$ many elements of $\LLL_{j,m,0}\cap [n',n'+p_m)$
we need to delete less than $2(d_{m}+1)\tqjm/q_{j',m}$
many for being too close to an element of $\LLL_{j',m,0}.$
We have a lower estimate
\begin{align}\label{*d5a}
\#(\LLL_{m}&\cap [n',n'+p_{m}))>
\sum_{j=1}^{K_{m}}\tqjm\bigg(1-\sum_{j'\not=j}2(d_{m}+1)\cdot \frac{1}
{q_{j',m}}\bigg)>\\ \nonumber
&\bigg (\sum_{j=1}^{K_{m}}\tqjm
\bigg )\bigg (1-K_{m}\frac{2(d_{m}+1)}
{\min_{j'} \{ q_{j',m}  \}}\bigg )=
p_{m}Q(m)\bigg (1-K_{m}\frac{2(d_{m}+1)}
{\min_{j'} \{ q_{j',m}  \}}\bigg )\\
& \nonumber >(1-\ggg_{m})p_{m}Q(m),
\end{align}
where,
taking into consideration \eqref{*f18},
the last inequality
for $m>3$
needs the assumption
\begin{equation}\label{*5aa}
\frac{1}{2000(m+1)
\oN_{m-2}}>K_{m}\frac{2(d_{m}+1)}{\min_{j'} \{ q_{j',m}  \}}
\end{equation}
about our initial parameters which can be achieved
by choosing the $q_{j,m}$'s sufficiently large.
For $m=2,3$ one needs to put $1/8$ to the left-hand side
of \eqref{*5aa} when this assumption is made.
Combining \eqref{*d4a} and \eqref{*d5a}
one can see that in any ``period" $[n',n'+p_{m})\sse [\bbb_{m-1},
\bbb_{m})$ the sequence $(n_{k})$ has a little less than
$p_{m}Q(m)$ many terms, and $Q(m)$ approximately
equals
the density
of this sequence here.
This can be reformulated as
\begin{equation}\label{**d5b}
1>\frac{\#(\LLL_{m}\cap [n',n'+p_{m}))}
{p_{m}Q(m)}
>1-\ggg_{m},
\end{equation}
or as
\begin{equation}\label{*d5b}
1>\frac{\#(\LLL_{m}\cap [n',n'+p_{m}))}
{\sum_{j=1}^{K_{m}}\#(\LLL_{j,m,0}\cap
[n',n'+p_{m}))}>1-\ggg_{m}.
\end{equation}

Later we need some assumptions and
estimations about our parameters.
In the rest of this section we give some of these,
not too difficult, but rather technical calculations.

We can choose our initial parameters
so that for all $m>0$ with $\ggg_{\bbb}$
defined in \eqref{*f18} we have
\begin{equation}\label{**f3c}
\bbb_{m-1}+2p_{m}<\frac{\ggg_{\bbb }}{2} \bbb_{m}.
\end{equation}
This implies
\begin{equation}\label{*f3c}
\bbb_{m}(1-\ggg_{\bbb })<(\bbb_{m}-\bbb_{m-1}-2p_{m}).
\end{equation}
Set $P_{m}=\lf \frac{\bbb_{m}-\bbb_{m-1}}{p_{m}}\rf$.
By \eqref{**d5b}
\begin{equation}\label{*f3d}
1>\frac{\#(\LLL_{m}\cap[\bbb_{m-1},\bbb_{m-1}+P_{m}p_{m}))}
{P_{m}p_{m}Q(m)}>1-\ggg_{m},
\end{equation}
and by \eqref{*d4a} we also have
\begin{equation}\label{*f3e}
\#(\LLL_{m}\cap [\bbb_{m-1},\bbb_{m-1}+P_{m}p_{m}))\leq
\oN_{\bbb_{m-1}}^{\bbb_{m}}<
\end{equation}
$$\#(\LLL_{m}\cap [\bbb_{m-1},\bbb_{m-1}+P_{m}p_{m}))+
p_{m}Q(m).$$

We need more estimates
of
$\oN_{\bbb_{m-1}}^{\bbb_{m}}$ from above, and
from below. By \eqref{*f3d} and \eqref{*f3e}
\begin{equation}\label{*f4aa}
\oN_{\bbb_{m-1}}^{\bbb_{m}}>(1-\ggg_{m})P_{m}p_{m}Q(m)=
\end{equation}
$$(1-\ggg_{m}) \bigg
\lf \frac{\bbb_{m}-\bbb_{m-1}}
{p_{m}} \bigg \rf p_{m}Q(m)>$$
(using \eqref{*f3c})
$$(1-\ggg_{m})((\bbb_{m}-\bbb_{m-1})-p_{m})Q(m)>
(1-\ggg_{m})(1-\ggg_{\bbb })\bbb_{m}Q(m),$$
on the other hand,
\begin{equation}\label{*f4ab}
\oN_{\bbb_{m-1}}^{\bbb_{m}}<P_{m}p_{m}Q(m)+p_{m}Q(m)
=(P_{m}+1)p_{m}Q(m)<
\end{equation}
(using \eqref{*f4c})
$$(\bbb_{m}-\bbb_{m-1}+p_{m})Q(m)<\bbb_{m}Q(m).$$

We suppose that an $m_{0}$ is given and
$\bbb_{m_{0}-1}< N \leq\bbb_{m_{0}}$.
Set $P_{m_{0},N}=\lf \frac{N-\bbb_{m_{0}-1}}
{p_{m_{0}}} \rf.$
By \eqref{**d5b}
\begin{equation}\label{*f5b}
1\geq\frac{\#(\LLL_{m_{0}}\cap
[\bbb_{m_{0}-1},\bbb_{m_{0}-1}+P_{m_{0},N}p_{m_{0}}))}
{P_{m_{0},N}p_{m_{0}}Q(m_{0})}>
1-\ggg_{m_{0}},
\end{equation}
where we regard $0/0=1$ by definition.
We also have
\begin{equation}\label{*f5a}
\#(\LLL_{m_{0}}\cap
[\bbb_{m_{0}-1},\bbb_{m_{0}-1}+P_{m_{0},N}p_{m_{0}}))
\leq \oN_{\bbb_{m_{0}}-1}^{N}<
\end{equation}
$$\#(\LLL_{m_{0}}\cap
[\bbb_{m_{0}-1},\bbb_{m_{0}-1}+P_{m_{0},N}p_{m_{0}}))
+p_{m_{0}}Q(m_{0}),$$
which implies
\begin{equation}\label{*f6aa}
\oN_{\bbb_{m_{0}-1}}^{N}
\geq(1-\ggg_{m_{0}})P_{m_{0},N}p_{m_{0}}Q(m_{0})
>
\end{equation}
$$(1-\ggg_{m_{0}})(N-\bbb_{m_{0}-1}-p_{m_0})Q(m_{0}),$$
and, on the other hand
\begin{equation}\label{*f6ab}
\oN_{\bbb_{m_{0}-1}}^{N}
<(P_{m_{0},N}+1)p_{m_{0}}Q(m_{0})\leq
(N-\bbb_{m_{0}-1}+p_{m_{0}})Q(m_{0}).
\end{equation}
To estimate $\oN_{0}^{N}$
from below
we combine \eqref{*f4aa} for $
m<m_{0}$ with $\eqref{*f6aa}$
\begin{equation}\label{*f4bb}
\oN_{0}^{N}=\sum_{m=1}^{m_{0}-1}\oN_{\bbb_{m-1}}^{\bbb_{m}}
+\oN_{\bbb_{m_{0}-1}}^{N}\geq
\end{equation}
$$\sum_{m=1}^{m_{0}-1}
(1-\ggg_{m})(\bbb_{m}-\bbb_{m-1}-p_{m})Q(m)
+(1-\ggg_{m_{0}})(N-\bbb_{m_{0}-1}-p_{m_{0}})Q(m_{0})>
$$
(using  \eqref{*f4c},
\eqref{*f6aa} and
$Q(m-1)\geq Q(m)$, $m=2,3,...$)
$$\frac{3}4\bigg (
\sum_{m=1}^{m_{0}-1}\frac{99}{100}(\bbb_{m}-\bbb_{m-1})Q(m)+
(N-\bbb_{m_{0}-1})Q(m_{0})-\frac{1}{100}(\bbb_{m_{0}-1}-\bbb_{m_{0}-2})
Q(m_{0})
\bigg )$$
$$>\frac{3}4Q(m_{0})\frac{98}{100}N>\frac{3}5Q(m_{0})N.$$

\section{Functions on $\Z$}\label{secfunz}

Assume $\fff:\Z\to\C$ is of finite support and $|\fff|\leq M$.
For ease of notation in this section we drop the subscript
$m$ corresponding to the $m$'th step of the definition
of $(n_{k})$. So we assume that $q_{1},...,q_{K}$ are
different primes and $p=q_{1}\cdot \cdot \cdot q_{K}.$
We will consider the $[(t-1)p,tp)\cap \Z$ ``grid intervals".
Given $n\in\Z$ choose $t(n)$ such that
$n\in[(t(n)-1)p,t(n)p).$
(In case there is a possibility of
 misunderstanding we will write
$t\cdot (n+2)$ for the product of $t$ and $(n+2)$
and $t(n+2)$ for the function $t(.)$ evaluated at $n+2$.)
For any $t\in\Z$ set
$\fff_{t,0}(n)=\fff(n-(t(n)-t)p)$. This function
is periodic by $p$ and coincides with $\fff$
on $[(t-1)p,tp)\cap \Z$, hence
\begin{equation}\label{*8cc}
\fff_{t(n),0}(n)=\fff(n).
\end{equation}
We also put
$$\offf_{0}(n)=\frac{1}p\sum_{k=(t(n)-1)p}^{t(n)p-1}
\fff(k)=\frac{1}p\sum_{k=0}^{p-1}\fff_{t(n),0}(k),$$
so $\offf_{0}(n)$ is the average of $\fff$ on the interval
$[(t(n)-1)p,t(n)p).$
Observe that
\begin{equation}\label{*d6}
\sum_{n=-\oo}^{\oo}\offf_{0}(n)=
\sum_{n=-\oo}^{\oo}\fff(n).
\end{equation}

Set
\begin{equation}\label{***d14}
t_{0}(n,N)=\bigg \lf \frac{n}p \bigg \rf+1=t(n),\ t_{1}(n,N)=
\bigg \lf
\frac{n+N}p \bigg \rf+1, \text{ and }
\end{equation}
$$N'=t_{1}(n,N)-t_{0}(n,N)+1.$$
For given $n$ and $N$ set
\begin{equation}\label{*d7}
I(n,N)=\bigg
[(t_{0}(n,N)-1)p-n,t_{1}(n,N)p-n\bigg )\cap \Z,
\end{equation}
\begin{equation}\label{**d7}
\nu(n,N)=\# I(n,N)=N'\cdot p,\ \nu(n,N,j)=\nu(n,N)/q_{j}=N'\tq_{j}.
\end{equation}
Clearly, $\nu(n,N)\leq N+p.$
We keep assumption \eqref{*d20f1}, that is,
\begin{equation}\label{*d7a}
1/2<q_{j}/q_{j'}<2,\text{ for any }j,j'.
\end{equation}
We introduce the operators
\begin{align}\label{*d9}
\cab(\fff,n,N,j)&=\frac1{\nu(n,N,j)}
\sum_{lq_{j}\in I(n,N)}\fff(n+lq_{j}),\\
\nonumber
\cab(\fff,n,N)&=
\frac{\sum_{j=1}^{K}\nu(n,N,j)\cab(\fff,n,N,j)}
{\sum_{j=1}^{K}\nu(n,N,j)},
\end{align}
and their ``modified versions"
\begin{align}\label{4*d7}
\cab_{0}(\fff,n,N,j)&=\frac1{\nu(n,N,j)}
\sum_{t=t_{0}(n,N)}^{t_{1}(n,N)}
\left |\sum_{lq_{j}+n\in [(t-1)p,tp)}\fff(n+lq_{j})
-\offf_{0}(n+lq_{j})\right|,\\
\nonumber
\cab_{0}(\fff,n,N)&=
\frac{\sum_{j=1}^{K}\nu(n,N,j)\cab_{0}(\fff,n,N,j)}
{\sum_{j=1}^{K}\nu(n,N,j)}.
\end{align}
Using (\ref{*d7}-\ref{*d7a}) it is not difficult to see
that
\begin{equation}\label{*d7b}
|\cab_{0}(\fff,n,N)|\leq\frac{2}K \sum_{j=1}^{K}
|\cab_{0}(\fff,n,N,j)|.
\end{equation}
The corresponding maximal operators are defined
as
$$\cab_{0}^{*}(\fff,n,j)=\sup_{N\geq 1}|\cab_{0}(\fff,n,N,j)|,
\text{ and }
\cab_{0}^{*}(\fff,n)=\sup_{N\geq 1}|\cab_{0}(\fff,n,N)|.$$

One of the main tools we will use later is the next lemma.
\begin{lemma}\label{d8}
For any $\fff:\Z\to\C$ of finite support, which is
bounded by $M$ we have
\begin{equation}\label{*ld8}
||\cab_{0}^{*}(\fff,.)||_{\ell ^{2}}\leq \frac{32}{K}
M ||\fff||_{\ell ^{1}}.
\end{equation}
\end{lemma}
The most useful
 ingredient in \eqref{*ld8} will be $K$ in the denominator
of the right-hand side.

In some estimates Lemma \ref{d8} will be used instead
of the usual maximal inequality (Lemma 3.5, p. 62 of \cite{[RW]}):
\begin{lemma}\label{Lemma3.5}
For all $\fff:\Z\to\C$ of finite support
for all $\lll>0$,
\begin{equation}\label{*d10}
\# \bigg \{ n:\sup_{N\geq 1}\left |\frac{1}N
\sum_{k=0}^{N-1}\fff(n+k)\right|>\lll \bigg  \}
\leq \frac{2}{\lll}||\fff||_{\ell^{1}}.
\end{equation}
\end{lemma}

We will also need the strong maximal inequality
from Lemma 4.4 of \cite{[RW]}. Here we give only the
special case of this lemma concerning $\ell ^{2}$
norms, and use slightly different notation.

\begin{lemma}\label{*strong}
For any $\fff:\Z\to\C$ of finite support
\begin{equation}\label{*eqstrong}
\bigg |\bigg |
\sup_{N>0}\frac{1}N\sum_{k=1}^{N}\fff(n+k)
\bigg |\bigg |_{\ell ^{2}}\leq 2||\fff||_{\ell ^{2}}.
\end{equation}
\end{lemma}

We recall a few basic facts about discrete Fourier
transforms.

For ease of notation we put $e(x)=\exp(2\pi i x).$

Given a function $\fff: \{ 0,...,p-1  \}\to \C$
we set
\begin{equation}\label{*FT1}
\hfff(\frac{b}p)=\frac{1}p\sum_{n=0}^{p-1}
\fff(n)e(-\frac{nb}p)\text{ for }b=0,...,p-1.
\end{equation}
Since $e(x)$ is periodic by one the definition
of $\hfff(b/p)$ can be extended for all
$b\in \Z$.

The inverse Fourier transform of
$\psi: \{ 0,\frac{1}p,...,\frac{p-1}p  \}\to\C$
is
\begin{equation}\label{*FT2}
\vppp(n)=\sum_{b=0}^{p-1}
\psi(\frac{b}p)e(n\frac{b}p)\text{ for }n=0,...,p-1.
\end{equation}
The way $\hfff$ and $\vppp$ are normalized differ in some
treatments, sometimes the factor
$1/p$ is used in the definition of $\vppp$
and sometimes factors of $1/\sqrt p$ are used in both definitions
of $\hfff$ and $\vppp$. With our choice of normalization
Parseval's theorem can be stated as
\begin{equation}\label{*FT3}
\frac{1}p\sum_{n=0}^{p-1}|\fff(n)|^{2}=
\sum_{b=0}^{p-1}|\hfff(\frac{b}p)|^{2}.
\end{equation}

Next we turn to the proof of Lemma \ref{d8}.

\begin{proof}[Proof of Lemma \ref{d8}.]
Set
\begin{equation}\label{*d11}
\fff_{t,0,j}(n)=\frac{1}{\tq_{j}}\sum_{k=0}^{\tq_{j}-1}
\fff_{t,0}(n+kq_{j}).
\end{equation}
This function is periodic by $q_{j}$ while
$\fff_{t,0}$ is periodic by $p=q_{j}\tq_{j}.$
The Fourier transform of $\fff_{t,0}$ is
$$\hfff_{t,0}(\frac{b}p)=\frac{1}p\sum_{n=0}^{p-1}\fff_{t,0}
(n)e(-\frac{nb}p),$$
while the Fourier transform of $\fff_{t,0,j}$ equals
\begin{equation}\label{*10aa}
 \hfff_{t,0,j}(\frac{b}p)=
\hfff_{t,0}(\frac{b}p)
\frac{1}{\tq_{j}}\sum_{k=0}^{\tq_{j}-1}
e(\frac{kq_{j}b}p)=
\left\{ \begin{array}{rll}
\hfff_{t,0}(\frac{b}p), & \mbox{if} & \tq_{j}|b; \\
0, & \mbox{if} & \tq_{j}\!\!\not|b.
\end{array}\right.
\end{equation}
Recall that $q_{j}$ and $q_{j'}$ are different primes when
$j\not=j'.$ Hence $0< b=r\tq_{j}=rp/q_{j}=r'\tq_{j'}=r'p/q_{j'}<p$
with integers $0<r<q_{j}$ and $0<r'<q_{j'}$ would imply
$rq_{j'}=r'q_{j}$, but this is impossible.
Since $\hfff_{t,0,j}$ is periodic by one
from this it follows  that for $b/p\not=0$
(modulo one)
and $j\not=j'$ we have
\begin{equation}\label{**d11}
\hfff_{t,0,j}(\frac{b}p)\hfff_{t,0,j'}(\frac{b}p)=0.
\end{equation}

Suppose $n\in [(t-1)p,tp)$.
Then
$\hfff_{t,0,j}(0)=\hfff_{t,0}(0)=\offf_{0}(n)$.
Set
$$\fff_{t,0,j,-}(n)=\fff_{t,0,j}(n)-\hfff_{t,0,j}(0)=
\fff_{t,0,j}(n)-\offf_{0}(n),$$
and
$$\fff_{t,0,-}(n)=\fff_{t,0}(n)-\hfff_{t,0}(0)=
\fff_{t,0}(n)-\offf_{0}(n)=\fff(n)-\offf_{0}(n)$$
where, again,
in the first display
the last equality and in the
second display the last two equalities  hold when
$n\in [(t-1)p,tp)$, that is, $t=t(n)$
while the other equalities make sense for other
$n$'s as well.

It is also clear that
\begin{equation}\label{**10aa}
\hfff_{t,0,-}(\frac{b}p)=\hfff_{t,0}(\frac{b}p)
\text{ if }b/p\not=0 \text{ mod }1, \text{ and }
\hfff_{t,0,-}(0)=0.
\end{equation}
We also put
$$\fff_{0,j,-}(n)=\fff_{t(n),0,j,-}(n),\text{ and }
\fff_{0,j,-}^{*}(n)=\sup_{N>0}\frac{1}N
\sum_{k=0}^{N-1}|\fff_{0,j,-}(n+kp)|.$$
By the strong maximal inequality
(Lemma \ref{*strong}) used on $n+p\Z$ instead of $\Z$,
$$\sum_{k=-\oo}^{\oo}|\fff_{0,j,-}^{*}(n+kp)|^{2}
\leq 2\sum_{k=-\oo}^{\oo}|\fff_{0,j,-}(n+kp)|^{2}.$$
Therefore,
\begin{equation}\label{*d12}
\sum_{n=0}^{p-1}\sum_{k=-\oo}^{\oo}
|\fff_{0,j,-}^{*}(n+kp)|^{2}=
\sum_{n=-\oo}^{\oo}|\fff_{0,j,-}^{*}(n)|^{2}
\leq
2\sum_{n=-\oo}^{\oo}|\fff_{0,j,-}(n)|^{2}.
\end{equation}
By Parseval's theorem and \eqref{**10aa}
$$\frac{1}p \sum_{n=0}^{p-1}
|\fff_{0,j,-}(n+(t-1)p)|^{2}=
\sum_{b=1}^{p}|\hfff_{t,0,j}(\frac{b}p)|^{2},$$
where we recall that $\hfff_{t,0,j,-}(0)=0$,
so this term is left out from the summation
on the right-hand side of the above formula.
It was the main motivation for introducing
the operators $\cab_{0}$, functions
$\fff_{t,0,j,-}$ and $\fff_{t,0,-}$.
Thus, keeping $t$ fixed
\begin{align}\label{*d13}
\frac{1}K\sum_{j=1}^{K}\sum_{n=0}^{p-1}
|\fff_{0,j,-}(n+(t-1)p)|^{2}&=
\frac{p}K\sum_{j=1}^{K}\sum_{b=1}^{p}
|\hfff_{t,0,j}(\frac{b}p)|^{2}\leq\\
\intertext{(using \eqref{*8cc},
(\ref{*d11}-\ref{**10aa}), and Parseval's theorem)}
\nonumber \frac{p}K\sum_{b=1}^{p}|\hfff_{t,0,-}(\frac{b}p)|^{2}
=&\frac{1}K\sum_{n=0}^{p-1}|\fff_{t,0,-}(n+(t-1)p)|^{2}=\\
\nonumber
\frac{1}K\sum_{n=0}^{p-1}|\fff(n+&(t-1)p)-\offf_{0}(n+(t-1)p)|^{2}.
\end{align}

Next we show that
\begin{equation}\label{**d13}
\fff^{*}_{0,j,-}(n)=\cab_{0}^{*}(\fff,n,j).
\end{equation}
By (\ref{***d14}-\ref{**d7}), $\nu(n,N)=N'p$, $\nu(n,N,j)=N'\tq_{j}$,
and by its definition
\begin{align}\label{*d14}
\cab_{0}(\fff,n,N,j)=
\frac{1}{\nu(n,N,j)}
\sum_{t=t_{0}(n,N)}^{t_{1}(n,N)}
\bigg |
\sum_{
lq_{j}+n\in
[(t-1)p,tp)
}
&\fff(n+lq_{j})-\offf_{0}(n+lq_{j})
\bigg |
=\\
\nonumber
\frac{1}{N'}\sum_{k=0}^{N'-1}
\bigg | \frac{1}{\tq_{j}}
\sum_{lq_{j}+n\in
[(t(n+kp)-1)p,t(n+kp)p)} (\fff(n+lq_{j})-&\offf_{0}(n+kp))\bigg |=
\\
\label{**d14}
\frac{1}{N'}\sum_{k=0}^{N'-1}|\fff_{t(n+kp),0,j}(n+kp)-
\offf_{0}(n+kp)|=&
\frac{1}{N'}\sum_{k=0}^{N'-1}|\fff_{0,j,-}(n+kp)|.
\end{align}
Taking supremum with respect to
$N$ in \eqref{*d14}, which means taking supremum with
respect to $N'$ in \eqref{**d14}, we obtain
\eqref{**d13}.

Clearly, by
(\ref{*d7}-\ref{*d7a}) and \eqref{4*d7}
$$|\cab_{0}^{*}(\fff,n)|\leq
\frac{\sum_{j=1}^{K}\nu(n,N,j)|\cab_{0}^{*}
(\fff,n,j)|}
{\sum_{j=1}^{K}\nu(n,N,j)}
\leq
\frac{2}K
\sum_{j=1}^{K}|\cab_{0}^{*}(\fff,n,j)|.$$
Using
this, \eqref{**d13} and the Cauchy-Schwarz
inequality
$$|\cab_{0}^{*}(\fff,n)|^{2}\leq
\frac{4}{K^{2}}\left (\sum_{j=1}^{K}
|\cab_{0}^{*}(\fff,n,j)|\right)^{2}\leq$$
$$\frac{4}{K^{2}}\left (\sqrt K\sqrt{\sum_{j=1}^{K}
|\cab_{0}^{*}(\fff,n,j)|^{2}}\right)^{2}\leq$$
$$\frac{4}K\sum_{j=1}^{K}|\cab_{0}^{*}(\fff,n,j)|^{2}
=\frac{4}{K}\sum_{j=1}^{K}|\fff_{0,j,-}^{*}(n)|^{2}.$$
Therefore, using \eqref{*d12}
$$\sum_{n=-\oo}^{\oo}|\cab_{0}^{*}(\fff,n)|^{2}
\leq \frac{4}K \sum_{j=1}^{K}\sum_{n=-\oo}^{\oo}
|\fff_{0,j,-}^{*}(n)|^{2}\leq
\frac{8}K\sum_{j=1}^{K}\sum_{n=-\oo}^{\oo}
|\fff_{0,j,-}(n)|^{2}=$$
$$\frac{8}K\sum_{j=1}^{K}\sum_{t=-\oo}^{\oo}
\sum_{n=0}^{p-1}|\fff_{0,j,-}(n+(t-1)p)|^{2}\leq$$
(using \eqref{*d13})
$$
\sum_{t=-\oo}^{\oo}\frac{8}K\sum_{n=0}^{p-1}
|\fff_{t,0,-}(n+(t-1)p)|^{2}=
\frac{8}K\sum_{n=-\oo}^{\oo}|\fff(n)-\offf_{0}(n)|^{2}\leq
$$
(using \eqref{*d6})
$$\frac{16}K\sum_{n=-\oo}^{\oo}(|\fff(n)|^{2}+|\offf_{0}(n)|^{2})
\leq \frac{16}K M\sum_{n=-\oo}^{\oo}(|\fff(n)|+|\offf_{0}(n)|)=
\frac{32}K M||\fff||_{\ell ^{1}}.$$
\end{proof}

We also need a weak $(1,1)$ inequality for the operator
$\cab^{*}$ which is defined as follows:
$$\cab^{*}(\fff,n)=\sup_{N>0}|\cab(\fff,n,N)|\text{ and we also use }
\cab^{*}(\fff,n,j)=\sup_{N>0}|\cab(\fff,n,N,j)|.$$

By using the definition of $\cab(\fff,n,N)$ and \eqref{*d7a}
it is easy to see that if $\fff\geq 0$ then
\begin{equation}\label{*d24}
\cab(\fff,n,N)\leq \frac{2}K\sum_{j=1}^{K}\cab(\fff,n,N,j).
\end{equation}

\begin{lemma}\label{dl24}
For any $\fff:\Z\to\R$
of finite support
and any $\tlll>0$ we have
\begin{equation}\label{**d24}
\# \{ n:\cab^{*}(\fff,n)>\tlll  \}\leq \frac{4}{\tlll}
||\fff||_{\ell ^{1}}.
\end{equation}
\end{lemma}
\begin{proof}
Since $\cab^{*}(|\fff|,n)\geq \cab^{*}(\fff,n)$ and
the right-hand side of
\eqref{**d24} is unchanged if $|\fff|$ is used instead of
$\fff$ we can suppose that $\fff\geq 0.$
Using notation introduced in the proof of Lemma \ref{d8}
set $\fff_{0,j,+}(n)=\fff_{t(n),0,j}(n).$
Then
using $\fff\geq 0$ one can see that for a fixed $t$
$$\sum_{n\in [(t-1)p,tp)}\fff(n)=
\sum_{n=(t-1)p}^{(t-1)p+q_{j}-1}
\sum_{k=0}^{\tq_j-1}\fff(n+kq_{j})=
$$
$$\sum_{n=(t-1)p}^{(t-1)p+q_{j}-1}\tq_j\fff_{t(n),0,j}(n)=
\sum_{n\in [(t-1)p,tp)}\fff_{0,j,+}(n),$$
and hence
\begin{equation}\label{*d17a}
||\fff_{0,j,+}||_{\ell ^{1}}=||\fff||_{\ell ^{1}}.
\end{equation}
Put
$$\cab_{K}^{*}(\fff,n)\defeq \sup_{N}\frac{1}K
\sum_{j=1}^{K}\cab(\fff,n,N,j),$$
$$\fff_{0,+}(n)=\frac{1}K\sum_{j=1}^{K}\fff_{0,j,+}(n),$$
and
$$\fff_{0,+}^{*}(n)=\sup_{N'}\frac{1}{N'}\sum_{k=0}^{N'-1}
\fff_{0,+}(n+kp).$$

Next we verify that
\begin{equation}\label{*d25}
\cab_{K}^{*}(\fff,n)=\fff_{0,+}^{*}(n).
\end{equation}

We use an argument similar to the one used at \eqref{**d13}.
Recall that
$$\cab_{K}^{*}(\fff,n)=
\sup_{N}\frac{1}K\sum_{j=1}^{K}
\frac{1}{\nu(n,N,j)}\sum_{  lq_{j}\in I(n,N)  }
\fff(n+lq_{j}).$$

Define $t_{0}(n,N),$ $t_{1}(n,N)$ and $N'$ as at \eqref{***d14}.
We have
\begin{align}\label{*d26}
\frac{1}{\nu(n,N,j)}&\sum_{  lq_{j}\in I(n,N)  }
\fff(n+lq_{j})=\\ \nonumber
&\frac{1}{N'\tq_{j}}\sum_{t=t_{0}(n,N)}^{t_{1}(n,N)}
\
\sum_{  lq_{j}+n\in[(t-1)\cdot p,t\cdot p)  }
\fff(n+lq_{j})=\\
\nonumber
&\frac{1}{N'}\sum_{k=0}^{N'-1}\frac{1}{\tq_{j}}
\sum_{  lq_{j}+n\in
[(t(n+kp)-1)\cdot p,t(n+kp)\cdot p ) }\fff(n+lq_{j})
=\\
\nonumber
&\frac{1}{N'}\sum_{k=0}^{N'-1}\fff_{t(n+kp),0,j}
(n+kp)=\frac{1}{N'}\sum_{k=0}^{N'-1}\fff_{0,j,+}
(n+kp).
\end{align}
Thus,
$$\frac{1}K\sum_{j=1}^{K}
\frac{1}{\nu(n,N,j)}
\sum_{  lq_{j}\in I(n,N)  }\fff(n+lq_{j})=
\frac{1}{N'}\sum_{k=0}^{N'-1}\frac{1}K\sum_{j=1}^{K}
\fff_{0,j,+}(n+kp).$$
Now taking supremums in $N$ and hence in $N'$ we obtain
\eqref{*d25}.

By \eqref{*d24} and \eqref{*d25}
$$\cab^{*}(\fff,n)\leq \sup_{N}
\frac{2}K \sum_{j=1}^{K}\cab(\fff,n,N,j)=$$
$$2\cab^{*}_{K}(\fff,n)=2\sup_{N'>0}\frac{1}{N'}
\sum_{k=0}^{N'-1}\frac{1}K\sum_{j=1}^{K}
\fff_{0,j,+}(n+kp)=2\fff_{0,+}^{*}(n).$$
Hence,
$$\# \{ n:\cab^{*}(\fff,n)>\tlll  \}\leq
\#  \{ n:\fff_{0,+}^{*}>\tlll/2  \}\leq$$
(using Lemma \ref{Lemma3.5} for $n+p\Z$ instead of $\Z$ and then adding for
$n$'s)
$$2\cdot \frac{2}{\tlll}
\sum_{n=0}^{p-1}\sum_{k\in\Z}\fff_{0,+}(n+kp)=
\frac{4}{\tlll}\sum_{n\in\Z}\fff_{0,+}(n)=$$
(using \eqref{*d17a})
$$\frac{4}{\tlll}\sum_{n\in\Z}\frac{1}K
\sum_{j=1}^{K}\fff_{0,j,+}(n)=\frac{4}{\tlll}
\frac{1}K\sum_{j=1}^{K}||\fff_{0,j,+}||_{\ell ^{1}}=
\frac{4}{\tlll}||\fff||_{\ell ^{1}}.$$
\end{proof}

In the next two sections we prove Lemma \ref{*dl2}.

\section{Part 1 of the proof of Lemma \ref{*dl2}}\label{secpart1}

\begin{proof}[Proof of Lemma \ref{*dl2}]
It is sufficient to show the lemma by assuming
$f\geq 0$ and by approximating $L^{1}$ functions
with simple functions we can assume that
$f$ takes finitely many values.

Suppose $\lll>0$ is fixed.
Set $X(A^{*})= \{ x: A^{*}(f,x)>\lll  \}$, and
$\lll'=\lll/3.$

Set
\[
f_{1,m}(x)=\left\{ \begin{array}{rll}
f(x)/\lll'& \mbox{if} & f(x)/\lll'<\oN_{m-3}; \\
0 &  \mbox{otherwise,} &
\end{array}\right.
\]
\[
f_{2,m}(x)=\left\{ \begin{array}{rll}
f(x)/\lll'& \mbox{if} &\oN_{m-3}\leq f(x)/\lll'<\oN_{m}; \\
0 &  \mbox{otherwise,} &
\end{array}\right.
\]
and
\[
f_{3,m}(x)=\left\{ \begin{array}{rll}
f(x)/\lll'& \mbox{if} &\oN_{m}\leq f(x)/\lll'; \\
0 &  \mbox{otherwise.} &
\end{array}\right.
\]
It is clear that
\begin{equation}\label{*d18b}
\sum_{m=1}^{\oo}f_{2,m}(x)\leq 3 f(x)/\lll',
\end{equation}
and
\begin{equation}\label{*15aa}
f_{1,m}=0 \text{ if }m\leq 3.
\end{equation}
If $x\in X(A^{*})$ then there exists $N$ such that
$A(f,x,N)>\lll.$
Since $3f(x)/\lll=
f(x)/\lll'=
f_{1,m}(x)+f_{2,m}(x)+f_{3,m}(x)$
if $m(N)$ is chosen so that
$\bbb_{m(N)-1}< N\leq\bbb_{m(N)}$
then there exists $i\in \{ 1,2,3  \}$
such that
$A(f_{i,m(N)},x,N)>1.$
Set $$X(A,m,i)= \{ x:\sup_{\bbb_{m-1}< N\leq\bbb_{m}} A(f_{i,m},x,N)>1 \},$$
and
$$X(A,m)= \{ x:\sup_{\bbb_{m-1}< N\leq\bbb_{m}} A(f,x,N)>\lll \}.$$
Then
\begin{equation}\label{*d18}
X(A^{*})=\cup_{m=1}^{\oo}X(A,m)\sse
\cup_{m=1}^{\oo}\cup_{i=1}^{3}X(A,m,i).
\end{equation}

The functions $n(x)$ and $r(x,m)$ will be defined later.
At this stage of the proof we only assume that they are
measurable in $x$, $r(.,m):X\to \{ 0,1,...,p_{m}-1\}$,
$n(.):X\to\Z.$
Given $N$ we let
\begin{equation}\label{*16bb}
I(x,m,N)=\bigg [-r(x,m),\bigg \lf \frac{N+r(x,m)}
{p_{m}} \bigg \rf p_{m}+p_{m}-r(x,m)\bigg )\cap \Z,
\end{equation}
$$\nu(x,m,N)=\# I(x,m,N),\qquad
\nu(x,m,N,j)=\nu(x,m,N)/q_{j,m}.$$
From \eqref{*d20f1} it follows that
\begin{equation}
\label{*d20f1b}
\frac{\nu(x,m,N,j)}
{\sum_{j'=1}^{K_{m}}
\nu(x,m,N,j')}
=\frac{\nu(x,m,N)\frac{1}{q_{j,m}}}
{\nu(x,m,N)Q(m)}\leq \frac{2}{K_{m}}.
\end{equation}

For any $g$ defined on $X$ we set
$$B(g,x,m,N,j)=
\frac{1}{\nu(x,m,N,j)}
\sum_{lq_{j,m}\in I(x,m,N)}g(T^{lq_{j,m}}x),$$
and
$$B(g,x,m,N)=
\frac{\sum_{j=1}^{K_{m}}
\nu(x,m,N,j)B(g,x,m,N,j)}
{\sum_{j=1}^{K_{m}}\nu(x,m,N,j)}.$$

From \eqref{*d20f1b} it follows that
for $g\geq 0$
\begin{equation}\label{*d20}
B(g,x,m,N)\leq \frac{2}{K_{m}}\sum_{j=1}^{K_{m}
}B(g,x,m,N,j).
\end{equation}

We also introduce the operator
$$A(g,x,N,m(N))=\frac{1}{\oN_{0}^{N}}
\sum_{k=\oN_{m(N)-1}+1}^{\oN_{0}^{N}}
g(T^{n_{k}}x),$$
and for $1\leq m<m(N)$ the operators
$$A(g,x,N,m)=\frac{1}{\oN_{0}^{N}}
\sum_{k=\oN_{m-1}+1}^{\oN_{m}}
g(T^{n_{k}}x).$$

We have
\begin{equation}\label{**d27}
A(g,x,N)=\sum_{m=1}^{m(N)}A(g,x,N,m).
\end{equation}

Next we verify that for any choice of $r(x,m)$,
$n(x)$, for any $g\geq 0$, $N\in\N$, if $m_{0}=m(N)$
then
\begin{equation}\label{*d27}
2 B(g,x,m_{0},N)\geq A(g,x,N,m_{0}).
\end{equation}

It is clear that $\nu(x,m_{0},N)\leq N+2p_{m_{0}}$ and by
$N>\bbb_{m_{0}-1},$ \eqref{*f4c}, and \eqref{*f4bb} we have
\begin{equation}\label{*f4b}
\frac{Q(m_{0})\nu(x,m_{0},N)}{\oN_{0}^{N}}<
\frac{Q(m_{0})1.01\cdot N}{\frac{3}5 Q(m_{0})N}<2.
\end{equation}

Now, still supposing $g\geq 0$
$$B(g,x,m_{0},N)=
\frac{\sum_{j=1}^{K_{m_{0}}}\sum_{lq_{j,m_{0}}\in I(x,m_{0},N)}
g(T^{lq_{j,m_{0}}}x)}
{\sum_{j=1}^{K_{m_{0}}}\nu(x,m_{0},N,j)}\geq$$
$$
\frac{\sum_{n_{k}\in [\bbb_{m_{0}-1},N)}g(T^{n_{k}}x)}
{Q(m_{0})\nu(x,m_{0},N)}\geq
A(g,x,N,m_{0})\cdot \frac{\oN_{0}^{N}}{Q(m_{0})\nu(x,m_{0},N)},
$$
that is,
$$
\frac{Q(m_{0})\nu(x,m_{0},N)}{\oN_{0}^{N}}B(g,x,m_{0},N)\geq
A(g,x,N,m_{0}),
$$
and \eqref{*f4b} implies \eqref{*d27}.

Suppose that for an $N\in(\bbb_{m_{0}-1},\bbb_{m_{0}}]$
we have $A(f_{2,m_{0}},x,N)>1$, that is,
$x\in X(A,m_{0},2).$

For $m\leq m_{0}$ and $\bbb_{m_{0}-1}<N'\leq N$ set
$$X(f_{2,m_{0}},N',m,\frac{1}3)=
 \{ x:A(f_{2,m_{0}},x,N',m)>\frac{1}3  \},$$
and
$$X(f_{2,m_{0}},N',m,+)=
 \{ x:A(f_{2,m_{0}},x,N',m)>0  \}.$$
Recall that if $f_{2,m_{0}}(x)\not=0$
then
$$\oN_{m_{0}-3}\leq f_{2,m_{0}}(x)<\oN_{m_{0}}.$$
We also put
$$X(f_{2,m_{0}},+)= \{ x:f_{2,m_{0}}(x)>0  \}.$$
Then
$$\cup_{m=1}^{m_{0}-3}X(f_{2,m_{0}},N,m,+)\sse
\cup_{k=1}^{\oN_{m_{0}-3}}T^{-n_{k}}(X(f_{2,m_{0}},+)),$$
which implies
\begin{equation}\label{*d28}
\mu(\cup_{m=1}^{{m_{0}-3}}X(f_{2,m_{0}},N,m,+))
\leq
\oN_{m_{0}-3}\mu(X(f_{2,m_{0}},+))\leq
\int f_{2,m_{0}}d\mu.
\end{equation}
For $m\leq m_{0}-1$ and any $N,N'\in (\bbb_{m_{0}-1},\bbb_{m_{0}}]$
$$X(f_{2,m_{0}},N,m,+)
=X(f_{2,m_{0}},N',m,+).$$
Thus from \eqref{*d28} we infer that
\begin{equation}\label{**d28}
\mu( \{ x:\sup_{\bbb_{m_{0}-1}< N\leq\bbb_{m_{0}}}
\sum_{m=1}^{m_{0}-3}A(f_{2,m_{0}},x,N,m)>0  \})
\leq \int f_{2,m_{0}}d\mu.
\end{equation}

Next we have to estimate
\begin{equation}\label{**d29}
\mu( \{ x:\sup_{\bbb_{m_{0}-1}< N\leq\bbb_{m_{0}}}
A(f_{2,m_{0}},x,N,m)>\frac{1}3  \})
\end{equation}
for $m=m_{0}-2,$ $m_{0}-1$ and $m_{0}.$

If $m'=m_{0}-2,$ or $m_{0}-1$
then for any $\bbb_{m_{0}-1}< N\leq \bbb_{m_{0}}$
we have $\bbb_{m'}<N$ and
$A(f_{2,m_{0}},x,N,m')\leq A(f_{2,m_{0}},x,\bbb_{m'},m')$.
Hence
\begin{equation}\label{*d29}
\mu( \{ x:\sup_{\bbb_{m_{0}-1}< N\leq\bbb_{m_{0}}}
A(f_{2,m_{0}},x,N,m')>\frac{1}3  \})
\leq
\end{equation}
$$\mu( \{ x:A(f_{2,m_{0}},x,\bbb_{m'},m')>\frac{1}3  \})\leq$$
$$3\int \frac{1}{\oN_{m'}}\sum_{k=\oN_{m'-1}}^{\oN_{m'}}
f_{2,m_{0}}(T^{n_{k}}x)d\mu(x)\leq
3\int f_{2,m_{0}}(x)d\mu(x).$$

To estimate \eqref{**d29} when $m=m_{0}$ is a little
more involved.

If $\int f_{2,m_{0}}d\mu=0$ then we have nothing to prove.
Hence, suppose
\begin{equation}\label{***d30}
\eee_{m_{0}}=\int f_{2,m_{0}}d\mu>0.
\end{equation}
Later we will choose a sufficiently large
$\kkk_{m_{0}}$ and by the Kakutani-Rokhlin lemma
a set $E_{2,m_{0}}$ such that
$E_{2,m_{0}},...,T^{\kkk_{m_{0}}-1}E_{2,m_{0}}$
are disjoint and
\begin{equation}\label{**d30}
\mu(\bigcup_{k=0}^{\kkk_{m_{0}}-1}T^{k}E_{2,m_{0}})
>1-\eee_{m_{0}}.
\end{equation}
Then $1/\mu(E_{2,m_{0}})\leq 1/\kkk_{m_{0}}$ and we can
assume that $\kkk_{m_{0}}$ is so large that
\begin{equation}\label{*d30}
(\bbb_{m_{0}}+3p_{m_{0}})\mu(E_{2,m_{0}})\leq
(\bbb_{m_{0}}+3p_{m_{0}})/\kkk_{m_{0}}<\eee_{m_{0}}.
\end{equation}
Since $f$ takes only finitely many values so does
$f_{2,m_{0}}$. Thus, we can partition each
$T^{k}E_{2,m_{0}}$ into a finite partition
$\aaa_{2,m_{0},k}$ so that $f_{2,m_{0}}$ is
constant on each partition element.
Consider $\aaa_{2,m_{0}}=\vee_{k=0}^{\kkk_{m_{0}}-1}
T^{-k}\aaa_{2,m_{0},k}.$
If $E'\in \aaa_{2,m_{0}}$
then $f_{2,m_{0}}$ is constant on each
$T^{k}E'$, $k=0,...,\kkk_{m_{0}}-1.$
It is enough to deal with the $E'$'s when
$\mu(E')>0$, and hence we suppose this.

Choose an arbitrary $x\in E'$ and set
$$\fff_{E'}(n)=f_{2,m_{0}}(T^{n}x)\text{ for }n=0,...,\kkk_{m_{0}}-1.$$
For other $n$'s set $\fff_{E'}(n)=0.$
If $x\not\in \cup_{n=0}^{\kkk_{m_{0}}-1}T^{n}E_{2,m_{0}},$
then set $r(x,m_0)=0.$
If $x\in \cup_{n=0}^{\kkk_{m_{0}}-1}T^{n}E_{2,m_{0}}$ then there is
a unique $E'(x)\in\aaa_{2,m_{0}}$ and $n(x)$ such that
$x\in T^{n(x)}E'(x)$.
In this case, set $r(x,m_0)=n(x)-\lf n(x)/p_{m_0} \rf p_{m_0}.$

Suppose $E'\in \aaa_{2,m_{0}}$ is fixed and $x\in \cup
_{n=0}^{\kkk_{m_{0}}-p_{m_{0}}-1}T^{n}E'$ and
$N\leq \bbb_{m_{0}}$.
Then letting $t(n(x))=\lf n(x)/p_{m_{0}} \rf+1$ we have
$n(x)\in [(t(n(x))-1)p_{m_{0}},t(n(x))p_{m_{0}})$
and if we use $p=p_{m_{0}}$ in \eqref{*d7} then
taking into consideration \eqref{*16bb}
$$I(n(x),N)=(t(n(x))-1)p_{m_{0}}
-n(x)+r(x,m_{0})+
I(x,m_{0},N)=I(x,m_{0},N),$$
$$\nu(n(x),N)=\nu(x,m_{0},N)\qquad \nu(n(x),N,j)=
\nu(x,m_{0},N,j),$$
and for $x\in \cup_{n=p_{m_{0}}}^{\kkk_{m_{0}}-\bbb_{m_{0}}-2p_{m_{0}}-1}
T^{n}E',$ with $q_{j}=q_{j,m_{0}}$  in the definition of $\cab$,
$$\cab(\fff_{E'},n(x),N,j)=B(f_{2,m_{0}},x,m_{0},N,j),$$
and
$$\cab(\fff_{E'},n(x),N)=B(f_{2,m_{0}},x,m_{0},N).$$

By \eqref{*d27} for $x\in \cup_{n=p_{m_{0}}}^{\kkk_{m_{0}}-\bbb_{m_{0}}-2p_{m_{0}}-1}
T^{n}E'$
\begin{equation}\label{*d32}
\sup_{\bbb_{m_{0}-1}< N\leq\bbb_{m_{0}}}A(f_{2,m_{0}},x,N,m_{0})\leq
2\sup_{\bbb_{m_{0}-1}< N\leq \bbb_{m_{0}}}B(f_{2,m_{0}},x,m_{0},N)=
\end{equation}
$$2\sup_{\bbb_{m_{0}-1}< N\leq\bbb_{m_{0}}}\cab(\fff_{E'},n(x),N)\leq
2\sup_{N>0}\cab(\fff_{E'},n(x),N)=2\cab^{*}(\fff_{E'},n(x)).$$
From \eqref{**d24} of Lemma \ref{dl24} it follows that
\begin{equation}\label{**d32}
\#\{n:\cab^{*}(\fff_{E'},n)>\frac{1}6\}\leq 24 \sum_{n\in\Z}\fff_{E'}(n).
\end{equation}

Using that $\mu(T^{n}E')=\mu(E')$ and the sets $T^{n}E'$
are disjoint for $n=0,...,\kkk_{m_{0}}-1$,
if we multiply both sides of \eqref{**d32}
by $\mu(E')$ and take into consideration \eqref{*d32}
then we obtain
\begin{equation}\label{*d33}
\mu
\left\{x\in\bigcup_{n=p_{m_{0}}}^{\kkk_{m_{0}}-\bbb_{m_{0}}-2p_{m_{0}}-1}
T^{n}E':\sup_{\bbb_{m_{0}-1}< N\leq\bbb_{m_{0}}}
A(f_{2,m_{0}},x,N,m_{0})>\frac{1}3
\right\}\leq\end{equation}
$$ 24 \int_{\cup_{n=0}^{\kkk_{m_{0}}
-1}T^{n}E'}f_{2,m_{0}}d\mu.$$
Adding \eqref{*d33} for all $E'\in\aaa_{2,m_{0}}$
we have
\begin{equation}\label{**d33}
\mu\bigg \{x\in\bigcup_{n=p_{m_{0}}}^{\kkk_{m_{0}}-\bbb_{m_{0}}-2p_{m_{0}}-1}
T^{n}E_{2,m_{0}}:\sup_{\bbb_{m_{0}-1}< N\leq\bbb_{m_{0}}}
A(f_{2,m_{0}},x,N,m_{0})>\frac{1}3\bigg \}\leq\end{equation}
$$ 24 \int_{X}f_{2,m_{0}}d\mu.$$
This, \eqref{***d30}, \eqref{**d30} and \eqref{*d30} imply that
\begin{equation}\label{***d33}
\mu\{x:\sup_{\bbb_{m_{0}-1}< N\leq \bbb_{m_{0}}}
A(f_{2,m_{0}},x,N,m_{0})>\frac{1}3\}\leq 26 \int f_{2,m_{0}}d\mu.
\end{equation}

Now,
\begin{equation}\label{*d34}
 \{ x:\sup_{N}A(f,x,N)>\lll  \}=
 \{ x:\exists N(x) ,\  A(f,x,N(x))>\lll \}=
\end{equation}
(we select and fix a measurable function $N(x)$)
$$= \{ x:\bbb_{m(N(x))-1}< N(x)\leq\bbb_{m_{N(x)}},\
A(f,x,N(x))>\lll  \}\sse$$
(by \eqref{*d18})
$$\bigcup_{m_{0}=1}^{\oo}\bigcup_{i=1}^{3}
 \{ x:\sup_{\bbb_{m_{0}-1}< N\leq \bbb_{m_{0}}}A(f_{i,m_{0}},x,N)>1  \}.$$
Using \eqref{**d27}
$$
 \{ x:\sup_{\bbb_{m_{0}-1}< N\leq \bbb_{m_{0}}}A(f_{2,m_{0}},x,N)>1  \}\sse
$$ $$
 \{ x:\sup_{\bbb_{m_{0}-1}< N\leq \bbb_{m_{0}}}
\sum_{m=1}^{m_{0}-3}A(f_{2,m_{0}},x,N,m)>0  \}\cup
$$
$$
\bigcup_{m'=m_{0}-2}^{m_{0}-1}
 \{ x:\sup_{\bbb_{m_{0}-1}< N\leq \bbb_{m_{0}}}A(f_{2,m_{0}},x,N,m')>\frac{1}3
  \}\cup
$$
$$
 \{ x:\sup_{\bbb_{m_{0}-1}< N\leq \bbb_{m_{0}}}A(f_{2,m_{0}},x,N,m_{0})
>\frac{1}3  \}.
 $$
This implies that
by \eqref{**d28}, \eqref{*d29}, and \eqref{***d33}
$$\mu \{ x:\sup_{\bbb_{m_{0}-1}< N\leq \bbb_{m_{0}}}A(f_{2,m_{0}},x,N)>1  \}
\leq 33\int f_{2,m_{0}}d\mu,$$
and
\begin{equation}\label{*d35}
\mu\left(\bigcup_{m_{0}=1}^{\oo}
 \{ x:\sup_{\bbb_{m_{0}-1}< N\leq \bbb_{m_{0}}}A(f_{2,m_{0}},x,N)>1  \}
\right )\leq 33 \int\sum_{m_{0}=1}^{\oo}f_{2,m_{0}}d\mu\leq
\end{equation}
(by \eqref{*d18b})
$$\frac{99\int fd\mu}{\lll'}\leq \frac{300\int fd\mu}{\lll}.$$

The estimation for the functions of the type $f_{3,m_{0}}$
is quite simple.
We have
\begin{align}\label{*d36}
\bigcup_{m_{0}=1}^{\oo}&
 \{ x:\sup_{\bbb_{m_{0}-1}< N\leq \bbb_{m_{0}}}A(f_{3,m_{0}},x,N)>1  \}
\sse\\ \nonumber &
\bigcup_{m_{0}=1}^{\oo}
 \{ x:\sup_{1\leq N\leq \bbb_{m_{0}}}A(f_{3,m_{0}},x,N)>0  \}
.
\end{align}
Put
\begin{equation}\label{***d36}
X_{3,m}= \{ x:\lll'\oN_{m+1}> f(x)\geq \lll'\oN_{m}  \}
.
\end{equation}
Observe that
$A(f_{3,m_{0}},x,N)=0$ if
$N\leq \bbb_{m_0}$ and
$x\not\in \cup_{k=1}^{\oN_{m_{0}}}
\cup_{m=m_{0}}^{\oo}T^{-n_{k}}X_{3,m}.$
Thus,
\begin{equation}\label{**d36}
\bigcup_{m_{0}=1}^{\oo}
 \{ x:\sup_{1\leq N\leq\bbb_{m_{0}}}A(f_{3,m_{0}},x,N)>0  \}\sse
\bigcup_{m_{0}=1}^{\oo}\bigcup_{k=1}^{\oN_{m_{0}}}\bigcup
_{m=m_{0}}^{\oo}T^{-n_{k}}X_{3,m}
=
\end{equation}
$$
\bigcup_{m_{0}=1}^{\oo}\bigcup
_{m=m_{0}}^{\oo}
\bigcup_{k=1}^{\oN_{m_{0}}}
T^{-n_{k}}X_{3,m}
=
\bigcup_{m_{0}=1}^{\oo}
\bigcup_{k=1}^{\oN_{m_{0}}}
T^{-n_{k}}X_{3,m_{0}}.
$$
From \eqref{*d36}, \eqref{***d36}, and \eqref{**d36}
it follows that
\begin{equation}\label{*d37}
\mu\left(\bigcup_{m_{0}=1}^{\oo}
 \{ x:\sup_{\bbb_{m_{0}-1}< N\leq \bbb_{m_{0}}}A(f_{3,m_{0}},x,N)>1  \}
\right )\leq
\mu \bigg (\bigcup_{m_{0}=1}^{\oo}
\bigcup_{k=1}^{\oN_{m_{0}}}
T^{-n_{k}}X_{3,m_{0}}\bigg )\leq
\end{equation}
$$\sum_{m_{0}=1}^{\oo}\oN_{m_{0}}\mu(X_{3,m_{0}})\leq
\int\frac{f(x)}
{\lll'}
d\mu(x)<\frac{3\int f d\mu}{\lll}.$$


\section{Part 2 of the proof of Lemma \ref{*dl2}}\label{secpart2}

Suppose $\bbb_{m_{0}-1}
< N \leq\bbb_{m_{0}}.$

We need to estimate
$\mu \{ x:\sup_{\bbb_{m_{0}-1}< N\leq \bbb_{m_{0}}}
A(f_{1,m_{0}},x,N)>1  \}.$
At the beginning we argue similarly to the case
$m=m_{0}$ when we had to obtain an estimate
of the functions $f_{2,m_{0}}$, however soon
this proof gets much more complicated. This is
mainly due to the fact that in the earlier argument
$\sum_{m_{0}}f_{2,m_{0}}$ is in $L^{1}$ while
we do not have this for $\sum_{m_{0}}f_{1,m_{0}}$.
To handle this problem
after we have applied
a Kakutani-Rokhlin tower construction
we need to take advantage
of the proper choice of constants $K_{m_{0}}$
and of Lemma \ref{d8}.

By choosing our initial parameters properly we can assume
that for all $m_{0}> 3,$
\begin{equation}\label{*d38f1}
(\sum_{m=1}^{m_{0}-2}\oN_{m})\frac{\oN_{m_{0}-3}}{\oN_{0}^{N}}
\leq
(\sum_{m=1}^{m_{0}-2}\oN_{m})\frac{\oN_{m_{0}-3}}{\oN_{m_{0}-1}}
<\frac{1}{100 m_{0}}.
\end{equation}

If $m_{0}\leq 3$ then $f_{1,m_{0}}=0$, hence
it is enough to obtain an estimate for $m_{0}> 3.$

By our assumptions and by its definition $0\leq
f_{1,m_{0}}<\oN_{m_{0}-3}$ and later
we will use this estimate quite often.

By \eqref{*d38f1} we have
\begin{equation}\label{*d38}
\sum_{m=1}^{m_{0}-2}A(f_{1,m_{0}},x,N,m)
\leq (\sum_{m=1}^{m_{0}-2}\oN_{m})\frac{\oN_{m_{0}-3}}
{\oN_{0}^{N}}<\frac{1}{100 m_{0}}.
\end{equation}

If $\int f_{1,m_{0}}d\mu=0$ then we have nothing to prove.
Hence, suppose
\begin{equation}\label{*d40}
\eee_{1,m_{0}}=\min\{2^{-m_{0}},2^{-m_{0}}\int f_{1,m_{0}}d\mu
\}>0.
\end{equation}
Later we will select a sufficiently large
$\kkk_{1,m_{0}}$
 and by the Kakutani-Rokhlin lemma
choose $E_{1,m_{0}}$ such that
$E_{1,m_{0}},...,T^{\kkk_{1,m_{0}}-1}E_{1,m_{0}}$
are disjoint and
\begin{equation}\label{**d40}
\mu(\cup_{k=0}^{\kkk_{1,m_{0}}-1}T^{k}E_{1,m_{0}})
>1-\eee_{1,m_{0}}.
\end{equation}
Then $1/\mu(E_{1,m_{0}})<1/\kkk_{1,m_{0}}$ and we can
assume that $\kkk_{1,m_{0}}$ is so large that
\begin{equation}\label{***d40}
(\bbb_{m_{0}}+3p_{m_{0}})\mu(E_{1,m_{0}})<
(\bbb_{m_{0}}+3p_{m_{0}})/\kkk_{1,m_{0}}<\eee_{1,m_{0}}.
\end{equation}
Since $f$ takes only finitely many values, so does
$f_{1,m_{0}}$. Thus we can divide each
$T^{k}E_{1,m_{0}}$ into a finite partition
$\aaa_{1,m_{0},k}$ so that $f_{1,m_{0}}$ is
constant on each partition element.
Consider $\aaa_{1,m_{0}}=\vee_{k=0}^{\kkk_{1,m_{0}}-1}
T^{-k}\aaa_{1,m_{0},k}.$
If $E'\in \aaa_{1,m_{0}}$
then $f_{1,m_{0}}$ is constant on each
$T^{k}E'$, $k=0,...,\kkk_{1,m_{0}}-1.$
It is enough to deal with the $E'$'s when
$\mu(E')>0$, and hence we suppose this.

Choose an arbitrary $x\in E'$ and set
$$\fff_{E'}(n)=f_{1,m_{0}}(T^{n}x)\text{ for }n=0,...,\kkk_{1,m_{0}}-1.$$
For other $n$'s set $\fff_{E'}(n)=0.$
If $x\not\in \cup_{n=0}^{\kkk_{1,m_{0}}-1}T^{n}E_{1,m_{0}},$
or $m>m_{0}$
then set $r(x,m)=0.$
If $x\in \cup_{n=0}^{\kkk_{1,m_{0}}-1}T^{n}E_{1,m_{0}}$ then there is
a unique $E'(x)\in\aaa_{1,m_{0}}$ and $n(x)$ such that
$x\in T^{n(x)}E'(x)$,
in this case
for $m'\leq m_{0}$
set $r(x,m')=n(x)-\lf n(x)/p_{m'} \rf p_{m'},$
$t(n(x),m')=\lf \frac{n(x)}{p_{m'}} \rf+1.$
This means that
$n(x)\in [(t(n(x),m')-1)p_{m'},t(n(x),m')p_{m'})$
and if we use $p=p_{m'}$ in \eqref{*d7} then by using \eqref{*16bb}
we have
$$I(n(x),N)= (t(n(x),m')-1)p_{m'}-n(x)+r(x,m')+
I(x,m',N)=I(x,m',N). $$
Still using $p=p_{m'}$ and $q_{j}=q_{j,m'}$
in \eqref{*d7} and \eqref{**d7} set
$$\nu(x,m',N)=\nu(n(x),N),\ \nu(x,m',N,j)=\nu(n(x),N,j).$$
We also put
$$t_{0}(x,m',N)=t(n(x),m'),\  t_{1}(x,m',N)=
t(n(x),m')+\nu(x,m',N)/p_{m'}.$$

For $x\in\cup_{n=p_{m'}}^{\kkk_{1,m_{0}}-p_{m'}-1}
T^{n}E'$ set
$$\off_{1,m_{0}}(x,m')=\frac{1}{p_{m'}}\sum_{k\in I(x,m',0)}
f_{1,m_{0}}(T^{k}x)=$$
$$\frac{1}{p_{m'}}\sum_{k=(t(n(x),m')-1)p_{m'}}
^{t(n(x),m')p_{m'}-1}\fff_{E'}(k)\defeq \offf_{E',m'}(n(x)).$$
If $x\not\in\cup_{n=p_{m'}}^{\kkk_{1,m_{0}}-p_{m'}-1}
T^{n}E'$ set $\off_{1,m_{0}}(x,m')=0.$

For $x\in\cup_{n=p_{m'}}^{\kkk_{1,m_{0}}-\bbb_{m_0}-p_{m_0}-1}
T^{n}E'$
and $0\leq N\leq \bbb_{m_0}$
set
$$B_{0}(f_{1,m_{0}},x,m',N,j)=
\frac{1}{\nu(x,m',N,j)}\cdot $$ $$\cdot
\sum_{t=t_{0}(x,m',N)}^{t_{1}(x,m',N)}
\bigg|\sum_{  lq_{j,m'}+n(x)\in [(t-1)p_{m'},tp_{m'})  }
f_{1,m_{0}}(T^{lq_{j,m'}}x)-\off_{1,m_{0}}(T^{lq_{j,m'}}x,m')\bigg|,$$
and
$B_{0}(f_{1,m_{0}},x,m',N)=\frac{\sum_{j=1}^{K_{m'}}
\nu(x,m',N,j)B_{0}(f_{1,m_{0}},x,m',N,j)}
{\sum_{j=1}^{K_{m'}}\nu(x,m',N,j)}.$
Observe that (for $N\leq\bbb_{m_{0}}$)
$$B_{0}(f_{1,m_{0}},x,m',N,j)=\cab_{0}
(\fff_{E'},n(x),N,j),$$
and hence
$$B_{0}(f_{1,m_{0}},x,m',N)=\cab_{0}(\fff_{E'},n(x),N),$$
provided $p=p_{m'}$, $q_{j}=q_{j,m'}$, $j=1,...,K_{m'}$
are used in the definition of $\cab_{0}$.
To emphasize this dependence on $m'$ we will use the notation
$$\cab_{0}(\fff_{E'},n(x),m',N,j)=
\cab_{0}(\fff_{E'},n(x),N,j),$$
and
$$\cab_{0}(\fff_{E'},n(x),m',N)=
\cab_{0}(\fff_{E'},n(x),N),$$
when the above choice of parameters is used.

 Set
\begin{align*}
&I_{1}(x,m_{0},N)=\Z\cap\\
&\bigg [\bigg \lf \frac{\bbb_{m_{0}-1}+r(x,m_{0})}
{p_{m_{0}}}\bigg \rf p_{m_{0}}-r(x,m_{0}),
\bigg \lf \frac{N+r(x,m_{0})}{p_{m_{0}}}\bigg \rf p_{m_{0}}+p_{m_{0}} -
r(x,m_{0})\bigg ),
\end{align*}
for $1\leq m<m_{0}$ set
\begin{align*}
&I_{1}(x,m,N)=\Z\cap\\
&\bigg [\bigg \lf \frac{\bbb_{m-1}+r(x,m)}
{p_{m}}\bigg \rf p_{m}-r(x,m),
\bigg \lf \frac{\bbb_{m}+r(x,m)}{p_{m}}\bigg \rf p_{m}+p_{m} -
r(x,m)\bigg ).
\end{align*}
We also put for $1\leq m\leq m_{0}$
$$\nu_{1}(x,m,N)=\# I_{1}(x,m,N),\qquad \nu_{1}(x,m,N,j)=
\frac{\nu_{1}(x,m,N)}{q_{j,m}}.$$

Next we need some estimates.
We also use the notation introduced in the end of
Section \ref{defnk}.
Clearly, for $m<m_{0}$
\begin{equation}\label{*f3cc}
P_{m}p_{m}\leq \nu_{1}(x,m,N)\leq (P_{m}+2)p_{m},
\end{equation}
and
\begin{equation}\label{*f3ccc}
P_{m_{0},N}p_{m_{0}}\leq \nu_{1}(x,m_{0},N)
\leq (P_{m_{0},N}+2)p_{m_{0}}.
\end{equation}
\begin{equation}\label{*f6ac}
N-\bbb_{m_{0}-1}\leq\nu_{1}(x,m_{0},N)\leq
(P_{m_{0},N}+2)p_{m_{0}}<N-\bbb_{m_{0}-1}+2p_{m_{0}}.
\end{equation}
By \eqref{*f4c} and \eqref{*f3c}
\begin{equation}\label{*f6ad}
\bbb_{m_{0}-1}(1-\ggg_{\bbb})<
\bbb_{m_{0}-1}-\bbb_{m_{0}-2}<
\nu_{1}(x,m_{0}-1,N)\leq
\end{equation}
$$(P_{m_{0}-1}+2)p_{m_{0}-1}<\bbb_{m_{0}-1}-\bbb_{m_{0}-2}+
2p_{m_{0}-1}<\bbb_{m_{0}-1},$$
From \eqref{*f6aa} and \eqref{*f6ab} it follows that
\begin{equation}\label{*f6ae}
|\oN_{\bbb_{m_{0}-1}}^{N}-Q(m_{0})(N-\bbb_{m_{0}-1})|<
\ggg_{m_{0}}Q(m_{0})(N-\bbb_{m_{0}-1})+p_{m_{0}}Q(m_{0}).
\end{equation}
On the other hand, by the definition of $I_{1}(x,m_{0},N)$
and $\nu_{1}(x,m_{0},N)$
\begin{equation}\label{*f6af}
|Q(m_{0})\nu_{1}(x,m_{0},N)-Q(m_{0})(N-\bbb_{m_{0}-1})|<
2p_{m_{0}}Q(m_{0}).
\end{equation}
Hence,
\begin{equation}\label{*f6b}
|\oN_{\bbb_{m_{0}-1}}^{N}-Q(m_{0})\nu_{1}(x,m_{0},N)|<
\ggg_{m_{0}}Q(m_{0})(N-\bbb_{m_{0}-1})+
3p_{m_{0}}Q(m_{0}).
\end{equation}
By \eqref{*f6aa} and \eqref{*f3ccc}
\begin{equation}\label{*f6ca}
\oN_{\bbb_{m_{0}-1}}^{N}>(1-\ggg_{m_{0}})P_{m_{0},N}p_{m_{0}}
Q(m_{0})>
\end{equation}
$$(1-\ggg_{m_{0}})(\nu_{1}(x,m_{0},N)-2p_{m_{0}})Q(m_{0}).$$
From \eqref{*f6ab} and \eqref{*f3ccc} it follows that
\begin{equation}\label{*f6cb}
\oN_{\bbb_{m_{0}-1}}^{N}<(\nu_{1}(x,m_{0},N)+p_{m_{0}})Q(m_{0}).
\end{equation}
Using \eqref{*f6ca} and \eqref{*f6cb} we infer
\begin{equation}\label{*f6cc}
|\oN_{\bbb_{m_{0}-1}}^{N}-\nu_{1}(x,m_{0},N)Q(m_{0})|<
\ggg_{m_{0}}\nu_{1}(x,m_{0},N)Q(m_{0})+
2p_{m_{0}}Q(m_{0}).
\end{equation}
By \eqref{*f4ab} and \eqref{*f6ad}
\begin{equation}\label{*f7a}
\oN_{\bbb_{m_{0}-2}}^{\bbb_{m_{0}-1}}<(\bbb_{m_{0}-1}-
\bbb_{m_{0}-2}+p_{m_{0}-1})Q(m_{0}-1)<
\end{equation}
$$(\nu_{1}(x,m_{0}-1,N)+p_{m_{0}-1})Q(m_{0}-1).$$
On the other hand, by \eqref{*f4aa} and \eqref{*f6ad}
\begin{equation}\label{*f7b}
\oN_{\bbb_{m_{0}-2}}^{\bbb_{m_{0}-1}}>(1-\ggg_{m_{0}-1})
(\nu_{1}(x,m_{0}-1,N)-3p_{m_{0}-1})Q(m_{0}-1).
\end{equation}
From \eqref{*f7a} and \eqref{*f7b} we infer
\begin{equation}\label{*f7c}
|\oN_{\bbb_{m_{0}-2}}^{\bbb_{m_{0}-1}}-\nu_{1}(x,m_{0}-1,N)
Q(m_{0}-1)|<
\end{equation}
$$\ggg_{m_{0}-1}\nu_{1}(x,m_{0}-1,N)Q(m_{0}-1)
+3p_{m_{0}-1}Q(m_{0}-1).$$

For $1\leq m\leq m_{0}$ set
$$S_{1}(f_{1,m_{0}},x,m,N,j)=\sum_{  lq_{j,m}\in I_{1}(x,m,N)
  } f_{1,m_{0}}(T^{lq_{j,m}}x),$$
$$\oS_{1}(f_{1,m_{0}},x,m,N,j)=\sum_{  lq_{j,m}\in I_{1}(x,m,N)
  } \off_{1,m_{0}}(T^{lq_{j,m}}x,m),$$
$$S_{1}(f_{1,m_{0}},x,m,N)=\sum_{j=1}^{K_{m}}
S_{1}(f_{1,m_{0}},x,m,N,j),$$
and
$$\oS_{1}(f_{1,m_{0}},x,m,N)=\sum_{j=1}^{K_{m}}
\oS_{1}(f_{1,m_{0}},x,m,N,j).$$

Until the end of the proof of this
lemma we assume that $m'=m_{0}-1$, or $m_{0}$.

Recall that in any subinterval of length $p_{m'}$
belonging to $[\bbb_{m'-1},\bbb_{m'})$ the sets
$\LLL_{j,m',0}$ have $\tq_{j,m'}=p_{m'}/q_{j,m'}$
many elements. From $\LLL_{j,m',0}$ during the definition
of $\LLL_{m'}$
(see \eqref{*d5a} and the paragraph above it)
less than
\begin{equation}\label{4*e6}
\sum_{j'\not=j}\tq_{j,m'}\frac{1}{q_{j',m'}}2(d_{m'}+1)
\end{equation}
many elements are deleted.
The intervals $[\bbb_{m_{0}-2},\bbb_{m_{0}-1})\cap\Z$
and $I_{1}(x,m_{0}-1,N)$ are roughly the same, apart
from two intervals of cardinality no more than $p_{m_{0}-1}$
at the beginning and in the end, to state this more
precisely
\begin{equation}\label{*e6}
\#(([\bbb_{m_{0}-2},\bbb_{m_{0}-1})\cap \Z)
\Delta I_{1}(x,m_{0}-1,N))\leq 2p_{m_{0}-1},
\end{equation}
where $\Delta$ stands for the symmetric difference.
Similarly,
\begin{equation}\label{*99bb}
\# (([\bbb_{m_{0}-1},N)\cap\Z)
\Delta I_{1}(x,m_{0},N))\leq 2p_{m_{0}},
\end{equation}
or, by changing by one element at the beginning and
in the end
\begin{equation}\label{**e6}
\#((\bbb_{m_{0}-1},N]\cap \Z)\Delta I_{1}(x,m_{0},N)<3 p_{m_{0}}.
\end{equation}
If $\can(m')$ denotes
 the total number of grid intervals of length $p_{m'}$
which are shifted by $-r(x,m')$ and are
belonging to
$I_{1}(x,m',N)$
then
\begin{equation}\label{***e6}
\can(m')=\frac{\nu_{1}(x,m',N,j)}{\tq_{j,m'}},
\text{ for any }j=1,...,K_{m'}.
\end{equation}

Next we verify that by our choice of the initial parameters we have
\begin{equation}\label{*d60a}
\sum_{j=1}^{K_{m'}}\nu_{1}(x,m',N,j)
<\sum_{j=1}^{K_{m'}}\nu(x,m',\min \{ N,\bbb_{m'}  \},j)
<2\oN_{0}^{N}.
\end{equation}
holds.

Observe that $\min \{ N,\bbb_{m'}  \}$ equals
$N$ when $m'=m_{0}$ and equals $\bbb_{m_{0}-1}$
if $m'=m_{0}-1$.

Since $N>\bbb_{m_{0}-1}$ and $m'\in  \{ m_{0},m_{0}-1  \}$
by \eqref{*f4c} we have
\begin{equation}\label{*xii27a}
N+p_{m'}<1.01\cdot N\text{ and }\bbb_{m_{0}-1}+p_{m'}<1.01\cdot
\bbb_{m_{0}-1}.
\end{equation}

By \eqref{**f3c} and \eqref{*f4bb} we have $$\oN_{0}^{N}\geq
\frac{3}4\cdot \frac{98}{100}(\bbb_{m_{0}-1}-\bbb_{m_{0}-2})
Q(m_{0}-1)>\frac{3}4\cdot \frac{98}{100}\cdot
\frac{999}{1000}\bbb_{m_{0}-1}Q(m_{0}-1)$$
and
$$\oN_{0}^{N}\geq \frac{3}5 Q(m_{0})N.$$

By the definition of $I(x,m',\min \{ N,\bbb_{m'}  \})$
and \eqref{*xii27a} we have
$$\nu(x,m',\min \{ N,\bbb_{m'}  \})<
\min \{ N,\bbb_{m'}  \}
+p_{m'}<1.01\cdot \min \{ N,\bbb_{m'}  \}.$$
Therefore
\begin{equation}\label{*xii27b}
\sum_{j=1}^{K_{m'}}\nu(x,m',\min \{ N,\bbb_{m'}  \}
,j)=\nu(x,m',\min \{ N,\bbb_{m'}  \})Q(m')
\leq 
\end{equation}
$$1.01\min \{ N,\bbb_{m'}  \}
Q(m')<2\oN_{0}^{N}.$$

We also make the following assumption about our
initial parameters:
\begin{equation}\label{*f15a}
\frac{\oN_{m_{0}-3}}{\oN_{m_{0}-1}}
3p_{m_{0}}
<\frac1{200m_{0}}\text{ for any }m_{0}> 3.
\end{equation}

From \eqref{*f15a} it follows that if $\bbb_{m_{0}-1}<N\leq \bbb_{m_{0}}$
then
\begin{equation}\label{*102b}
\frac{\oN_{m_{0}-3}}{\oN_{0}^{N}}3p_{m_{0}-1}<
\frac{\oN_{m_{0}-3}}{\oN_{m_{0}-1}}3p_{m_{0}}<\frac{1}{200 m_{0}}.
\end{equation}

Using \eqref{*f15} and (\ref{4*e6}-\ref{*102b})
for $m'=m_{0}-1,$ or $m_{0}$
we have
\begin{equation}\label{*d60}
\bigg |\frac{S_{1}(f_{1,m_{0}},x,m',N)}
{\oN_{0}^{N}}-A(f_{1,m_{0}},x,N,m')\bigg |\leq
\end{equation}
$$\frac{\oN_{m_{0}-3}}{\oN_{0}^{N}}
\bigg |3p_{m'}+\sum_{j=1}^{K_{m'}}
\nu_{1}(x,m',N,j)\sum_{j'\not=j}\frac{2(d_{m'}+1)}
{\min \{ q_{j',m'}  \}}\bigg |\leq$$
$$\frac{\oN_{m_{0}-3}}{\oN_{0}^{N}}
\bigg |3p_{m'}+\sum_{j=1}^{K_{m'}}
\nu_{1}(x,m',N,j)K_{m'}\frac{2(d_{m'}+1)}
{\min \{ q_{j',m'}  \}}\bigg |<\frac{1}{100m_{0}}.$$

Next we estimate
\begin{equation}\label{*d62}
\bigg |\frac{S_{1}(f_{1,m_{0}},x,m',N)-\oS_{1}(f_{1,m_{0}},x,m',N)}
{\oN_{0}^{N}}\bigg |=
\end{equation}
$$\bigg |\frac{\sum_{j=1}^{K_{m'}}\sum_{  lq_{j,m'}\in
I_{1}(x,m',N)}f_{1,m_{0}}(T^{lq_{j,m'}}x)-\off_{1,m_{0}}
(T^{lq_{j,m'}}x,m')  }
{\oN_{0}^{N}}\bigg |\leq$$
(using \eqref{*xii27b}, $I_{1}(x,m',N)\sse I(x,m',N)$ and the triangle
inequality)
\begin{align*}
&\frac{2}{\sum_{j=1}^{K_{m'}}\nu(x,m',\min \{ N,\bbb_{m'}\},j  )}\cdot
\\
&\cdot \sum_{j=1}^{K_{m'}}\sum_{t=t_{0}(x,m',\min \{ N,\bbb_{m'}  \})}
^{t_{1}(x,m',\min \{ N,\bbb_{m'}  \})}
\bigg |\sum_{  lq_{j,m'}+n(x)\in [(t-1)p_{m'},tp_{m'})  }
f_{1,m_{0}}(T^{lq_{j,m'}}x)-\\
&
\off_{1,m_{0}}(T^{lq_{j,m'}}x,m')\bigg |=
\end{align*}
$$2B_{0}(f_{1,m_{0}},x,m',\min \{ N,\bbb_{m'}  \})\leq
2\max_{\bbb_{m'-1}< N'\leq\bbb_{m'}}
B_{0}(f_{1,m_{0}},x,m',N').$$

Since $\off_{1,m_{0}}(T^{lq_{j,m'}}x,m')$
equals $\off_{1,m_{0}}(T^{t'p_{m'}}x,m')$
when $lq_{j,m'}\in [t'p_{m'},(t'+1)p_{m'})$
and for each
$t'$ and
$j$ there are $\tq_{j,m'}$
many such $lq_{j,m'}$'s we have
\begin{equation}\label{***d63}
\oS_{1}(f_{1,m_{0}},x,m',N)=
\sum_{j=1}^{K_{m'}}\tq_{j,m'}
\sum_{t'p_{m'}\in I_{1}(x,m',N)}\off_{1,m_{0}}(T^{t'p_{m'}}x,m')=
\end{equation}
$$\sum_{j=1}^{K_{m'}}\frac{p_{m'}}{q_{j,m'}}
\bigg (\frac{1}{p_{m'}}
\sum_{k\in I_{1}(x,m',N)}f_{1,m_{0}}(T^{k}x)\bigg )=
Q(m')\sum_{k\in I_{1}(x,m',N)}f_{1,m_{0}}(T^{k}x).$$
This implies
\begin{equation}\label{*106bb}
|\oS_{1}(f_{1,m_{0}},x,m',N)|<
\oN_{m_{0}-3}Q(m')\nu_{1}(x,m',N).
\end{equation}
We also have
\begin{equation}\label{*29bb}
\oN_{0}^{N}=\bigg (\sum_{m=1}^{m_{0}-2}\oN_{\bbb_{m-1}}^{\bbb_{m}}
\bigg )+\oN_{\bbb_{m_{0}-2}}^{\bbb_{m_{0}-1}}+
\oN_{\bbb_{m_{0}-1}}^{N},
\end{equation}
and
the initial parameters can be chosen so that
 we can estimate the sum on the right-hand side by
\begin{equation}\label{*d63}
\sum_{m=1}^{m_{0}-2}
\oN_{\bbb_{m-1}}^{\bbb_{m}}
< \frac{1}{100 m_{0}}
\oN_{\bbb_{m_{0}-2}}^{\bbb_{m_{0}-1}}.
\end{equation}

By \eqref{*f6ad}
and a suitable assumption about our initial parameters

\begin{equation}\label{*f8}
\nu_{1}(x,m_{0}-1,N)>(1-\ggg_{\bbb})\bbb_{m_{0}-1}
\text{ and }
\end{equation}
$$\frac{3p_{m_{0}-1}}
{\nu_{1}(x,m_{0}-1,N)}<\frac{3p_{m_{0}-1}}
{(1-\ggg_{\bbb})\bbb_{m_{0}-1}}<
\frac{1}{200 m_{0}\oN_{m_{0}-3}}.$$

From \eqref{*f18},
\eqref{*f7c} and \eqref{*f8} it follows that
\begin{equation}\label{**d63}
\bigg |\frac{1}{
\oN_{\bbb_{m_{0}-2}}^{\bbb_{m_{0}-1}}
}-\frac{1}{\nu_{1}(x,m_{0}-1,N)Q(m_{0}-1)}\bigg |=
\end{equation}
$$
\bigg |\frac{\nu_{1}(x,m_{0}-1,N)Q(m_{0}-1)-
\oN_{\bbb_{m_{0}-2}}^{\bbb_{m_{0}-1}}}
{\nu_{1}(x,m_{0}-1,N)Q(m_{0}-1)
\oN_{\bbb_{m_{0}-2}}^{\bbb_{m_{0}-1}}}\bigg |<
$$
$$\frac{1}{
\oN_{\bbb_{m_{0}-2}}^{\bbb_{m_{0}-1}}}
\bigg (\ggg_{m_{0}-1}+\frac{3p_{m_{0}-1}}
{\nu_{1}(x,m_{0}-1,N)}\bigg )
<\frac{1}{
\oN_{\bbb_{m_{0}-2}}^{\bbb_{m_{0}-1}}}
\cdot \frac{1}{100\cdot m_{0}\oN_{m_{0}-3}}.$$

Using \eqref{***d63} and \eqref{**d63}
\begin{equation}\label{**d64}
\bigg |\frac{\oS_{1}(f_{1,m_{0}},x,m_{0}-1,N)}
{\oN_{\bbb_{m_{0}-2}}^{\bbb_{m_{0}-1}}}
-\frac{1}{\nu_{1}(x,m_{0}-1,N)}
\sum_{k\in I_{1}(x,m_{0}-1,N)}f_{1,m_{0}}(T^{k}x)\bigg |<
\end{equation}
$$\oS_{1}(f_{1,m_{0}},x,m_{0}-1,N)
\frac{1}{\oN_{m_{0}-3}100m_{0}
\oN_{\bbb_{m_{0}-2}}^{\bbb_{m_{0}-1}}}<
$$
(By \eqref{*f4c}, \eqref{*f3c}, \eqref{*f7b} and \eqref{*106bb})
$$
\frac{\oN_{m_{0}-3}Q(m_{0}-1)\nu_{1}(x,m_{0}-1,N)}
{\oN_{m_{0}-3}100m_{0}\oN_{\bbb_{m_{0}-2}}^{\bbb_{m_{0}-1}}}
<\frac{2}{100 m_{0}}.$$

To obtain an estimate similar to
\eqref{**d64} for $m_{0}$ instead of $m_{0}-1.$
we separate two cases.\\
 CASE 1 holds if $N-\bbb_{m_{0}-1}\geq
10^{4}(m_{0}+1)\oN_{m_{0}-2}p_{m_{0}},$
and\\
 CASE 2 holds when
$0\leq N-\bbb_{m_{0}-1}<
10^{4}(m_{0}+1)\oN_{m_{0}-2}p_{m_{0}}.$

If CASE 1 holds by \eqref{***d63} we have
\begin{equation}\label{*f9}
\bigg | \frac{\oS_{1}(f_{1,m_{0}},x,m_{0},N)}
{\oN_{\bbb_{m_{0}-1}}^{N}}
-\frac{\sum_{k\in I_{1}(x,m_{0},N)}f_{1,m_{0}}(T^{k}x)}
{\nu_{1}(x,m_{0},N)}\bigg |=
\end{equation}
$$\oS_{1}(f_{1,m_{0}},x,m_{0},N)\bigg |\frac{1}{
\oN_{\bbb_{m_{0}-1}}^{N}}
-\frac{1}{\nu_{1}(x,m_{0},N)Q(m_{0})}\bigg |\leq$$
(using \eqref{*106bb})
$$Q(m_{0})\nu_{1}(x,m_{0},N)\oN_{m_{0}-3}
\frac{
|\oN_{\bbb_{m_{0}-1}}^{N}
-Q(m_{0})\nu_{1}(x,m_{0},N)|}
{
\oN_{\bbb_{m_{0}-1}}^{N}
Q(m_{0})\nu_{1}(x,m_{0},N)}=
$$ $$\frac{\oN_{m_{0}-3}|
\oN_{\bbb_{m_{0}-1}}^{N}
-Q(m_{0})\nu_{1}(x,m_{0},N)|}
{
\oN_{\bbb_{m_{0}-1}}^{N}}<$$
(using \eqref{*f6aa} and \eqref{*f6b})
$$\frac{\oN_{m_{0}-3}
(\ggg_{m_{0}}Q(m_{0})(N-\bbb_{m_{0}-1})
+3p_{m_{0}}Q(m_{0}))}
{(1-\ggg_{m_{0}})(N-\bbb_{m_{0}-1}-p_{m_{0}})Q(m_{0})}<$$
(using
\eqref{*f18} and that for CASE 1 we have
$\ggg_{m_{0}}(N-\bbb_{m_{0}-1})>3p_{m_{0}}$,
$N-\bbb_{m_{0}-1}>2p_{m_{0}}$ and $\ggg_{m_{0}}<1/2$)
$$\frac{\oN_{m_{0}-3}4\ggg_{m_{0}}}
{1-\ggg_{m_{0}}}<8\oN_{m_{0}-3}\ggg_{m_{0}}<
\frac{1}{100m_{0}}.$$

If CASE 2 holds then
\begin{equation}\label{*f12}
|A(f_{1,m_{0}},x,N,m_{0})|=\frac{1}{\oN_{0}^{N}}
\sum_{  n_{k}\in[\bbb_{m_{0}-1},N)  }
f_{1,m_{0}}(T^{n_{k}}x)<
\end{equation}
$$ \frac{1}{\oN_{0}^{N}}(N-\bbb_{m_{0}-1})\oN_{m_{0}-3}<
\frac{10^{4}(m_{0}+1)\oN_{m_{0}-2}
p_{m_{0}}
\oN_{m_{0}-3}
}
{\oN_{m_{0}-1}}<\frac{1}{1000{m_{0}}},$$
where the last inequality holds if a suitable
assumption is made about our initial parameters.

For both CASEs we also have
\begin{equation}\label{*d65}
\bigg |\frac{1}{\nu_{1}(x,m_{0}-1,N)}
\sum_{k\in I_{1}(x,m_{0}-1,N)}f_{1,m_{0}}(T^{k}x)-
\frac{1}{\bbb_{m_{0}-1}}\sum_{k=1}^{\bbb_{m_{0}-1}}
f_{1,m_{0}}(T^{k}x)\bigg |\leq
\end{equation}
$$\bigg |\frac{1}{\nu_{1}(x,m_{0}-1,N)}-
\frac{1}{\bbb_{m_{0}-1}}\bigg |\cdot
\sum_{k\in I_{1}(x,m_{0}-1,N)}
f_{1,m_{0}}(T^{k}x)+
$$ $$
\frac{1}{\bbb_{m_{0}-1}}\bigg |
\sum_{k=1}^{\bbb_{m_{0}-1}}f_{1,m_{0}}(T^{k}x)-
\sum_{k\in I_{1}(x,m_{0}-1,N)}f_{1,m_{0}}(T^{k}x)
\bigg |\leq$$
(using \eqref{*e6})
$$\frac{|\nu_{1}(x,m_{0}-1,N)
-\bbb_{m_{0}-1}|}
{\nu_{1}(x,m_{0}-1,N)
\bbb_{m_{0}-1}}\sum_{k\in I_{1}(x,m_{0}-1,N)}
f_{1,m_{0}}(T^{k}x)+
\frac{(\bbb_{m_{0}-2}+2p_{m_{0}-1})\oN_{m_{0}-3}}{\bbb_{m_{0}-1}}
\leq$$
(using that \eqref{*e6} implies $|\nu_{1}(x,m_{0}-1,N)-
\bbb_{m_{0}-1}|\leq \bbb_{m_{0}-2}+2p_{m_{0}-1}$)
$$\frac{2(\bbb_{m_{0}-2}+2p_{m_{0}-1})\oN_{m_{0}-3}
}{\bbb_{m_{0}-1}}<\frac{1}{100 m_{0}},$$
where at the last inequality we again made an assumption
about our initial parameters, especially we used that
$p_{m_{0}-1}<\bbb_{m_{0}-2}$ can be supposed to be much less than
$\bbb_{m_{0}-1}.$

Next observe that by \eqref{**e6}
\begin{equation}\label{*e11}
\bigg |\sum_{k\in I_{1}(x,m_{0},N)}f_{1,m_{0}}(T^{k}x)-
\bigg (\sum_{k=1}^{N}f_{1,m_{0}}(T^{k}x)-
\sum_{k=1}^{\bbb_{m_{0}-1}}f_{1,m_{0}}(T^{k}x)\bigg )\bigg |<
3p_{m_{0}}\oN_{m_{0}-3}.
\end{equation}
It is also clear
from \eqref{*99bb}
that
\begin{equation}\label{**e11}
|\nu_{1}(x,m_{0},N)
-(N-\bbb_{m_{0}-1})|\leq 2 p_{m_{0}},
\end{equation}
furthermore $p_{m_{0}}>1,$
$N-\bbb_{m_{0}-1}\geq 1 $ and \eqref{*f6ac} imply
\begin{equation}\label{***e11}
\nu_{1}(x,m_{0},N)
<3p_{m_{0}}(N-\bbb_{m_{0}-1}).
\end{equation}

By \eqref{**e11}
\begin{equation}\label{*115bb}
\bigg |\frac{1}{\nu_{1}(x,m_{0},N)}-
\frac{1}{N-\bbb_{m_{0}-1}}\bigg |=
\frac{|\nu_{1}(x,m_{0},N)
-(N-\bbb_{m_{0}-1})|}
{\nu_{1}(x,m_{0},N)
(N-\bbb_{m_{0}-1})}\leq
\end{equation}
$$\frac{2p_{m_{0}}}
{\nu_{1}(x,m_{0},N)
(N-\bbb_{m_{0}-1})}.$$

Hence,
\begin{equation}\label{*e12}
\bigg |\frac{1}{\nu_{1}(x,m_{0},N)}
\sum_{k\in I_{1}(x,m_{0},N)}f_{1,m_{0}}(T^{k}x)
-
\end{equation}
$$\frac{\sum_{k=1}^{N}f_{1,m_{0}}(T^{k}x)-\sum_{k=1}^{\bbb_{m_{0}-1}}
f_{1,m_{0}}(T^{k}x)}{N-\bbb_{m_{0}-1}}
\bigg |<
$$
(by using \eqref{*e11})
$$\frac{1}{\nu_{1}(x,m_{0},N)}3p_{m_{0}}\oN_{m_{0}-3}+
\bigg |\frac{1}{\nu_{1}(x,m_{0},N)}-\frac{1}{N-\bbb_{m_{0}-1}} \bigg |
(N-\bbb_{m_{0}-1})\oN_{m_{0}-3}\leq $$
(by \eqref{*f6ac}, \eqref{**e11} and \eqref{*115bb})
$$\frac{3p_{m_{0}}\oN_{m_{0}-3}}{(N-\bbb_{m_{0}-1})}
+\frac{2p_{m_{0}}}{\nu_{1}(x,m_{0},N)}\oN_{m_{0}-3}\leq$$
$$\frac{5 p_{m_{0}}\oN_{m_{0}-3}}
{N-\bbb_{m_{0}-1}}\defeq\cae.$$

If CASE 1 holds, that is, $N-\bbb_{m_{0}-1}\geq 10^{4}p_{m_{0}}
\oN_{m_{0}-2}(m_{0}+1)$ then
\begin{equation}\label{**e12}
\cae< \frac{1}{100m_{0}}.
\end{equation}

Otherwise, if CASE 2 holds then
\begin{equation}\label{***d67}
0<\oN_{0}^{N}-\oN_{\bbb_{m_{0}-2}}^{\bbb_{m_{0}-1}}<
\bbb_{m_{0}-2}
+10^{4}p_{m_{0}}\oN_{m_{0}-2}(m_{0}+1).
\end{equation}

By \eqref{*f7c}
\begin{equation}\label{*suppl0}
\bigg | \frac{\oN_{\beta _{m_0-2}}^{\beta _{m_{0}-1}}}
{Q(m_{0}-1)\nu_{1}(x,m_{0}-1,N)}-1\bigg |
< \ggg_{m_{0}-1}+ \frac{3p_{m_{0}-1}}{\nu_{1}(x,m_{0}-1,N)}<
\end{equation}
(using \eqref{*f18}, \eqref{*f4c} and \eqref{*f6ad})
$$<\frac{1}{2000}+\frac{3p_{m_{0}-1}}
{\bbb_{m_{0}-1}-\bbb_{m_{0}-2}}<\frac{1}{1000}.$$
Hence,
\begin{equation}\label{*suppl1}
\frac{Q(m_{0}-1)\nu_{1}(x,m_{0}-1,N)}
{\oN_{\bbb_{m_{0}-2}}^{\bbb_{m_{0}-1}}}<2.
\end{equation}

By \eqref{*106bb}
$$|\oS_{1}(f_{1,m_{0}},x,m_{0}-1,N)|<
\oN_{m_{0}-3}Q(m_{0}-1)\nu_{1}(x,m_{0}-1,N), $$
therefore,
\begin{equation}\label{*suppl2}
\bigg | \frac{\oS_{1}(f_{1,m_{0}},x,m_{0}-1,N)}
{\oN_{0}^{N}}-
\frac{\oS_{1}(f_{1,m_{0}},x,m_{0}-1,N)}
{\oN_{\bbb_{m_{0}-2}}^{\bbb_{m_{0}-1}}}
\bigg |<\end{equation}
$$\oN_{m_{0}-3}Q(m_{0}-1)\nu_{1}(x,m_{0}-1,N)
\frac{|\oN_{0}^{N}-\oN_{\bbb_{m_{0}-2}}^{\bbb_{m_{0}-1}}|}
{\oN_{0}^{N}\oN_{\bbb_{m_{0}-2}}^{\bbb_{m_{0}-1}}}<$$
(by \eqref{***d67} and \eqref{*suppl1})
$$\frac{\oN_{m_{0}-3}2(\bbb_{m_{0}-2}+10^{4}p_{m_{0}}
\oN_{m_{0}-2}(m_{0}+1))}
{\oN_{0}^{N}}\leq$$
$$\frac{\oN_{m_{0}-3}2(\bbb_{m_{0}-2}+10^{4}p_{m_{0}}
\oN_{m_{0}-2}(m_{0}+1))}
{\oN_{m_{0}-1}}<\frac{1}{200 m_{0}}$$
if a suitable assumption is made about our initial parameters.

Furthermore,
\begin{equation}\label{*d67}
\bigg |\frac{1}N\sum_{k=1}^{N}f_{1,m_{0}}(T^{k}x)-
\frac{1}{\bbb_{m_{0}-1}}\sum_{k=1}^{\bbb_{m_{0}-1}}
f_{1,m_{0}}(T^{k}x)\bigg |<
\end{equation}
$$\bigg |\frac{1}N-\frac{1}{\bbb_{m_{0}-1}}\bigg |
\sum_{k=1}^{\bbb_{m_{0}-1}}f_{1,m_{0}}(T^{k}x)+
\frac{1}N\sum_{k=\bbb_{m_{0}-1}+1}^{N}f_{1,m_{0}}(T^{k}x)<$$
$$\frac{N-\bbb_{m_{0}-1}}{N}
\left (\frac{1}{\bbb_{m_{0}-1}}\sum_{k=1}^{\bbb_{m_{0}-1}}
f_{1,m_{0}}(T^{k}x)\right)+\frac{1}N\oN_{m_{0}-3}(N-\bbb_{m_{0}-1})\leq
$$
(recalling that CASE 2 holds)
$$2\cdot \frac{N-\bbb_{m_{0}-1}}{N}
\oN_{m_{0}-3} <2\cdot 10^{4}p_{m_{0}}\oN_{m_{0}-2}(m_{0}+1)\frac{1}N
\oN_{m_{0}-3}<
$$
$$
2\cdot 10^{4} p_{m_{0}}\oN_{m_{0}-2}(m_{0}+1)\frac{1}{\bbb_{m_{0}-1}}
\oN_{m_{0}-3}
<\frac{1}{200 m_{0}},$$
if proper assumptions are made about our initial parameters.

To make easier to follow estimate \eqref{*d68}
in an abbreviated form we recall that\\ by
\eqref{*d60}, $|S_{1}/\oN_{0}^{N}-A|< 1/(100 m_{0})$,\\
by \eqref{*d62}, $|(S_{1}-\oS_{1})/\oN_{0}^{N}|\leq
2 \max B_{0}$,\\ by \eqref{*suppl2}, $|
(\oS_{1}/\oN_{0}^{N})-(\oS_{1}/\oN_{\bbb_{m_{0}-2}}
^{\bbb_{m_{0}-1}})|<1/(200 m_{0})$,\\
by \eqref{**d64}, $|(\oS_{1}/\oN_{\bbb_{m_{0}-2}}
^{\bbb_{m_{0}-1}})-(1/\nu_{1})\sum_{I_{1}} f_{1,m_{0}}|<
2/(100 m_{0})$,\\ by \eqref{*d65}, $|((1/\nu_{1})\sum_{I_{1}}
f_{1,m_{0}})-((1/\bbb_{m_{0}-1})\sum_{1}^{\bbb_{m_{0}-1}}f_{1,m_{0}})|<
1/(100 m_{0})$\\
 and by \eqref{*d67},
$|((1/N)\sum_{1}^{N}f_{1,m_{0}})-
((1/\bbb_{m_{0}-1})\sum_{1}^{\bbb_{m_{0}-1}}f_{1,m_{0}})|
<
1/(200 m_{0}).$

Thus in CASE 2 by \eqref{*d38}, \eqref{*d60}, \eqref{*d62},
\eqref{*29bb},
\eqref{**d64},
\eqref{*f12}, \eqref{*d65},
\eqref{*suppl2} and \eqref{*d67}
\begin{equation}\label{*d68}
\bigg |A(f_{1,m_{0}},x,N)-\frac{1}N\sum_{k=1}^{N}f_{1,m_{0}}(T^{k}x)\bigg |<
\bigg |\sum_{m=1}^{m_{0}-2}A(f_{1,m_{0}},x,N,m)\bigg |
+
\end{equation}
$$\bigg |A(f_{1,m_{0}},x,N,m_{0}-1)-\frac{1}N\sum_{k=1}^{N}
f_{1,m_{0}}(T^{k}x)\bigg |+|A(f_{1,m_{0}},x,N,m_{0})|<$$
$$
\frac{6}{100 m_{0}}+2\max_{\bbb_{m_{0}-2}< N'\leq\bbb_{m_{0}-1}}
B_{0}(f_{1,m_{0}},x,m_{0}-1,N')
+
\frac1{1000m_{0}}<$$ $$
\frac{1}{10m_{0}}
+
2\max_{\bbb_{m_{0}-2}< N'\leq\bbb_{m_{0}-1}}
B_{0}(f_{1,m_{0}},x,m_{0}-1,N').
$$

Next we need similar type estimates for CASE 1.

By
the assumption for
CASE 1, $N-\bbb_{m_{0}-1}\geq 10^{4}(m_{0}+1)
\oN_{m_{0}-2}p_{m_{0}}$,
moreover by \eqref{*f6ab},
$\oN_{\bbb_{m_{0}-1}}^{N}<(N-\bbb_{m_{0}-1}+p_{m_{0}})Q(m_{0})$,
and
by \eqref{*f4bb},
$\oN_{0}^{N}>\frac{3}5 N Q(m_{0}).$
Thus
\begin{equation}\label{**f13}
\frac{\oN_{\bbb_{m_{0}-1}}^{N}}{\oN_{0}^{N}}<
\frac{5}{3}
\frac{N-\bbb_{m_{0}-1}+p_{m_{0}}}
{N}<
\frac{5}{3}\frac{N-\bbb_{m_{0}-1}}{N}
(1+\frac{p_{m_{0}}}{N-\bbb_{m_{0}-1}})
\leq
\end{equation}
$$
 \frac{5}{3}
\frac{N-\bbb_{m_{0}-1}}
{N}(1+\frac{1}{
10^{4}(m_{0}+1)\oN_{m_{0}-2}})<
\frac{2(N-\bbb_{m_{0}-1})}
{N}<2.$$

If CASE 1 holds using
\eqref{*d60}, \eqref{*d62}, \eqref{*29bb}, \eqref{**d64}
and
\eqref{*d65} (see the list of abbreviated estimates before
\eqref{*d68} as well)
\begin{equation}\label{*d69}
\bigg |A(f_{1,m_{0}},x,N,m_{0}-1)-
\frac{\oN_{\bbb_{m_{0}-2}}^{\bbb_{m_{0}-1}}}
{\oN_{0}^{N}}\frac{1}{\bbb_{m_{0}-1}}\sum_{k=1}^{\bbb_{m_{0}-1}}
f_{1,m_{0}}(T^{k}x)\bigg |\leq
\end{equation}
$$2\max_{\bbb_{m_{0}-2}< N'\leq\bbb_{m_{0}-1}}
B_{0}(f_{1,m_{0}},x,m_{0}-1,N')
+\frac{1}{10 m_{0}}.$$

In addition to the list of abbreviated estimates
given before \eqref{*d68} we also recall
that\\ by
\eqref{*f9} we have
$|(\oS_{1}/\oN_{\bbb_{m_{0}-1}}^{N})-(1/\nu_{1})
(\sum_{k\in I_{1}}f_{1,m_{0}})|<1/(100 m_{0})$,\\
moreover by
\eqref{*e12} and \eqref{**e12} we have\\
$|(1/\nu_{1})(\sum_{k\in I_{1}}f_{1,m_{0}})-
(\sum_{k=1}^{N}f_{1,m_{0}}-\sum_{k=1}^{\bbb_{m_{0}-1}}f_{1,m_{0}})
/(N-\bbb_{m_{0}-1})|<1/(100 m_{0}).$

By \eqref{*d60}, \eqref{*d62}, \eqref{*29bb}, \eqref{*f9},
\eqref{*e12},
\eqref{**e12} and \eqref{**f13}
\begin{equation}\label{*e15}
\bigg |A(f_{1,m_{0}},x,N,m_{0})-\frac{\sum_{k=\bbb_{m_{0}-1}+1}^{N}
f_{1,m_{0}}(T^{k}x)}
{N-\bbb_{m_{0}-1}}\cdot \frac{\oN_{\bbb_{m_{0}-1}}^{N}}
{\oN_{0}^{N}}\bigg |<
\end{equation}
$$2\max_{\bbb_{m_{0}-1}< N'\leq\bbb_{m_{0}}}
B_{0}(f_{1,m_{0}},x,m_{0},N')
+\frac{1}{10 m_{0}}.$$
Set
$$X(f_{1,m_{0}},B_{0},m')=
\bigg  \{ x:\max_{\bbb_{m'-1}< N\leq\bbb_{m'}}
B_{0}(f_{1,m_{0}},x,m',N)
>\frac{1}{100\cdot 2^{m_{0}}} \bigg \}.$$
For $x\in\cup _{n=p_{m_{0}}}^{\kkk_{1,m_{0}}-\bbb_{m_{0}}-
p_{m_{0}}-1} T^{n}E'$
we have
\begin{equation}\label{*e16}
\max_{\bbb_{m'-1}< N\leq\bbb_{m'}}
B_{0}(f_{1,m_{0}},x,m',N)
= \max_{\bbb_{m'-1}< N\leq\bbb_{m'}}
\cab_{0}(\fff_{E'},n(x),m',N)
\leq
\end{equation}
$$ \sup_{0<N}
\cab_{0}(\fff_{E'},n(x),m',N)
=\cab_{0}^{*}(\fff_{E'},n(x),m').$$

By Lemma \ref{d8}
$$||\cab_{0}^{*}(\fff_{E'},.,m')||_{\ell ^{2}}\leq
\frac{32}{K_{m'}}\oN_{m_{0}-3}||\fff_{E'}||_{\ell ^{1}}.$$

Hence, (using $m'=m_{0}-1$, or $m_{0}$)
\begin{equation}\label{*d45}
\#
\bigg
\{ n:\cab_{0}^{*}(\fff_{E'},n,m')>\frac{1}{100\cdot 2^{m_{0}}}
\bigg
\}
\leq (100\cdot 2^{m_{0}})^{2}||\cab_{0}^{*}(\fff_{E'},.,m')||_{\ell ^{2}}\leq
\end{equation}
(using \eqref{*f13} for $m'=m_{0}-1$, or $m_{0}$)
$$10^{4}4^{m_{0}}\frac{32}{K_{m'}}\oN_{m_{0}-3}||\fff_{E'}||_{\ell ^{1}}
<2^{-m_{0}}\sum_{n\in\Z}\fff_{E'}(n).$$

Recalling that $\mu(T^{n}E')=\mu(E')$ and the sets
$T^{n}E'$ are disjoint for $n=0,...,\kkk_{1,m_{0}}-1$
if we multiply both sides of \eqref{*d45} by $\mu(E')$,
take into consideration that $\fff_{E'}(n)=0$
when $n\not\in  \{ 0,...,\kappa_{1,m_{0}-1}  \}$ and
we also use
 \eqref{*e16} we obtain
\begin{equation}\label{*d46}
\mu \bigg  \{
x\in\bigcup_{n=p_{m_{0}}}^{\kkk_{1,m_{0}}-\bbb_{m_{0}}-p_{m_{0}}-1}
T^{n}E':\max_{\bbb_{m'-1}< N\leq \bbb_{m'}}
B_{0}(f_{1,m_{0}},x,m',N)>\frac{1}{100\cdot 2^{m_{0}}} \bigg  \}
\leq
\end{equation}
$$2^{-m_{0}}\int_{\cup_{n=0}^{\kkk_{1,m_{0}}-1}T^{n}E'}
f_{1,m_{0}}d\mu.$$

Adding \eqref{*d46} for all $E'\in\aaa_{1,m_{0}}$ we have
\begin{equation}\label{**d46}
\mu \bigg
 \{ x\in\bigcup_{n=p_{m_{0}}}^{\kkk_{1,m_{0}}-\bbb_{m_{0}}-p_{m_{0}}-1}
T^{n}E_{1,m_{0}}:\max_{\bbb_{m'-1}< N\leq \bbb_{m'}}
B_{0}(f_{1,m_{0}},x,m',N)>\frac{1}{100\cdot 2^{m_{0}}} \bigg  \}
\leq
\end{equation}
$$2^{-m_{0}}\int
f_{1,m_{0}}d\mu.$$
This \eqref{*d40}, \eqref{**d40} and \eqref{***d40}
imply
\begin{equation}\label{*e18}
\mu(X(f_{1,m_{0}},B_{0},m'))=
\mu \bigg  \{ x:\max_{\bbb_{m'-1}< N\leq\bbb_{m'}}
B_{0}(f_{1,m_{0}},x,m',N)>\frac{1}{100\cdot 2^{m_{0}}} \bigg  \}
\leq\end{equation}
$$ 4\cdot 2^{-m_{0}}\int f_{1,m_{0}}d\mu.$$

Set $X(f,B_{0})=\cup_{m_{0}=1}^{\oo}(X(f_{1,m_{0}},B_{0},
m_{0}-1)\cup X(f_{1,m_{0}},B_{0},m_{0})).$
By \eqref{*e18}
\begin{equation}\label{**e18}
\mu(X(f,B_{0}))\leq {8}\sum_{m_{0}=1}^{\oo}2^{-m_{0}}
\int f_{1,m_{0}}d\mu\leq
\end{equation}
$$\frac{8}{\lll'}\int fd\mu={24}\frac{\int fd\mu}
{\lll}.$$

We also put
$$X(f,B_{0},\oo)=
\bigcap_{m=1}^{\oo}\bigcup_{m_{0}=m}^{\oo}
(X(f_{1,m_{0}},B_{0},m_{0}-1)
\cup X(f_{1,m_{0}},B_{0},m_{0})).$$
From \eqref{*e18} it follows that
\begin{equation}\label{***e18}
\mu(X(f,B_{0},\oo))=0.
\end{equation}

By the Wiener-Yosida-Kakutani Maximal Ergodic Theorem
if we set
$$X^{*}(f)= \bigg \{ x:\sup_{0<N}\frac{1}N
\sum_{k=1}^{N}f(T^{k}x)>\frac{\lll'}{100} \bigg  \}$$
then
\begin{equation}\label{*e19a}
\mu(X^{*}(f))<\frac{100}{\lll'}\int fd\mu=\frac{300}{\lll}
\int f d\mu.
\end{equation}

Suppose $x\in  X\sm(X^{*}(f)\cup X(f,B_{0}))$ and $N>0.$
Then there exists $m_{0}$ such that $\bbb_{m_{0}-1}< N
\leq \bbb_{m_{0}}$.

Since $f_{1,m_{0}}=0$ for $m_{0}\leq 3$ we can assume
$m_{0}>3.$

If CASE 2 holds then using \eqref{*d68} and $0\leq f_{1,m_{0}}\leq
f/\lll'$ we have
\begin{equation}\label{*e19}
A(f_{1,m_{0}},x,N)\leq \frac{1}{10m_{0}}+\frac{2}{100\cdot 2^{m_{0}}}
+\frac{1}{100}<1.
\end{equation}

If CASE 1 holds  for $x\in X\sm X^{*}(f)$
using $\oN_{\bbb_{m_{0}-2}}^{\bbb_{m_{0}-1}}
\leq \oN_{0}^{N}$
we have
$$\frac{\oN_{\bbb_{m_{0}-2}}^{\bbb_{m_{0}-1}}
}{\oN_{0}^{N}}\frac{1}{\bbb_{m_{0}-1}}\sum_{k=1}^{\bbb_{m_{0}-1}}
f_{1,m_{0}}(T^{k}x)<\frac{1}{100},$$
and hence by \eqref{*d69}
\begin{equation}\label{*e20}
A(f_{1,m_{0}},x,N,m_{0}-1)\leq
\frac{2}{100\cdot 2^{m_{0}}}+\frac{1}{10 m_{0}}+\frac{1}{100}
\end{equation}
for $x\in X\sm(X^{*}(f)\cup X(f,B_{0})).$

By
$f_{1,m_{0}}\geq 0$ and
\eqref{**f13} for $x\not\in X^{*}(f)$
\begin{equation}\label{**e20}
\bigg |\frac{\sum_{k=\bbb_{m_{0}-1}+1}^{N}
f_{1,m_{0}}(T^{k}x)}{N-\bbb_{m_{0}-1}}
\cdot \frac{\oN_{\bbb_{m_{0}-1}}^{N}}{\oN_{0}^{N}}\bigg |
\leq \frac{2}{N}\sum_{k=1}^{N}f_{1,m_{0}}(T^{k}x)\leq
\frac{2}{100}.
\end{equation}

Using \eqref{*e15} and \eqref{**e20} we obtain for
$x\in X\sm(X^{*}(f)\cup X(f,B_{0}))$
\begin{equation}\label{***e20}
A(f_{1,m_{0}},x,N,m_{0})<\frac{2}{100\cdot 2^{m_{0}}}+
\frac{1}{10 m_{0}}+\frac{2}{100}.
\end{equation}
From \eqref{*d38}, \eqref{*e20}, and \eqref{***e20} we infer
\begin{equation}\label{*e21a}
A(f_{1,m_{0}},x,N)\leq
\bigg (\sum_{m=1}^{m_{0}-2} A(f_{1,m_{0}},x,N,m)\bigg )
+ A(f_{1,m_{0}},x,N,m_{0}-1)+
\end{equation}
$$
A(f_{1,m_{0}},x,N,m_{0})<
\frac{1}{100m_{0}}+2\bigg (\frac{2}{100\cdot 2^{m_{0}}}
+\frac{1}{10 m_{0}}+\frac{2}{100}\bigg)<1.$$
Hence if $x\in X\sm(X^{*}(f)\cup X(f,B_{0}))$
for both CASEs
by \eqref{*e19}, or by \eqref{*e21a} we have
$$\sup_{\bbb_{m_{0}-1}< N\leq\bbb_{m_{0}}}A(f_{1,m_{0}},x,N)
<1$$
for any
$N\geq 1$ and $m_{0}$ satisfying $\bbb_{m_{0}-1}<
N\leq\bbb_{m_{0}}$,
and therefore by \eqref{**e18} and \eqref{*e19a}
\begin{equation}\label{*e21}
\mu\bigg(\bigcup_{m_{0}=1}^{\oo} \bigg \{ x:
\sup_{\bbb_{m_{0}-1}< N\leq\bbb_{m_{0}}}A(f_{1,m_{0}},x,N)
> 1 \bigg \}\bigg )\leq
\end{equation}
$$
\mu(X^{*}(f)\cup X(f,B_{0}))\leq
(300+{24})\frac{\int fd\mu}{\lll}.$$

Now \eqref{*d18}, \eqref{*d35}, \eqref{*d37} and \eqref{*e21}
imply
$$\mu \{ x:\sup_{0<N}A(f,x,N)>\lll  \}\leq
1000\frac{\int fd\mu}{\lll}.$$
This proves Lemma \ref{*dl2}.

\end{proof}

\section{The proof of Lemma \ref{dl1}}\label{sec7}

\begin{proof}[Proof of Lemma \ref{dl1}.]
We will use in this proof notation introduced
in the proof of Lemma \ref{*dl2}.
Without limiting generality we can assume $0\leq f\leq 1.$
To prove Lemma \ref{dl1} set $\lll=3$, that is,
$\lll'=1$ in the previous proof. Suppose
$N\geq \bbb_{4}$. Using $m_{0}=m(N)$, ($\bbb_{m_{0}-1}
\leq N<\bbb_{m_{0}}$) we have
$f_{1,m_{0}}(x)=f(x).$ Assume
$x\not\in X(f,B_{0},\oo).$
Then there exists $N(x,0,\oo)$ such that for $m_{0}\geq N(x,0,\oo)$,
$x\not\in X(f_{1,m_{0}},B_{0},m_{0}-1)\cup X(f_{1,m_{0}},B_{0},m_{0})=
X(f,B_{0},m_{0}-1)\cup X(f,B_{0},m_{0}).$
By the Ergodic Theoreom there exists $X^{**}(f)$ such that
$\mu(X^{**}(f))=0$ and if $x\not\in X^{**}(f)$ then
$\frac{1}{N}\sum_{k=1}^{N}f(T^{k}x)\to\int fd\mu.$

Suppose $\eee>0.$ If $x\not\in X^{**}(f)$ then there exists
$N(x,\eee)$ such that for $N\geq N(x,\eee)$ we have
$$\bigg |\frac{1}N\sum_{k=1}^{N}f(T^{k}x)-\int fd\mu\bigg |<\eee.$$
Suppose $x\not\in X(f,B_{0},\oo)\cup X^{**}(f)$ and
$N\geq N^{*}(x,\eee)=\max \{ N(x,0,\oo),N(x,\eee)  \}.$

If CASE 2 holds with $m_{0}=m(N)$ we obtain from
\eqref{*d68} that
$$|A(f,x,N)-\frac{1}{N}\sum_{k=1}^{N}
f(T^{k}x)|<\frac{1}{10m(N)}+\frac{2}{100\cdot 2^{m(N)}},$$
and hence
\begin{equation}\label{*e24}
|A(f,x,N)-
\int f d\mu|<
\eee+
\frac{1}{10m(N)}+\frac{2}{100\cdot 2^{m(N)}}.
\end{equation}

If CASE 1 holds with $m_{0}=m(N)$ we obtain from
\eqref{*d69}
$$\bigg |A(f,x,N,m(N)-1)-\frac{
\oN_{\bbb_{m(N)-2}}^{\bbb_{m(N)-1}}}{\oN_{0}^{N}}
\frac{1}{\bbb_{m(N)-1}}\sum_{k=1}^{\bbb_{m(N)-1}}f(T^{k}x)
\bigg |\leq \frac{2}{100\cdot 2^{m(N)}}+\frac{1}{10 m(N)},$$
which implies
\begin{equation}\label{**e24}
\bigg |A(f,x,N,m(N)-1)-\frac{
\oN_{\bbb_{m(N)-2}}^{\bbb_{m(N)-1}}}{\oN_{0}^{N}}
\int fd\mu
\bigg |\leq
\end{equation}
$$ \frac{2}{100\cdot 2^{m(N)}}+\frac{1}{10 m(N)}+
\frac{\oN_{\bbb_{m(N)-2}}^{\bbb_{m(N)-1}}}
{\oN_{0}^{N}}\eee
< \frac{1}{50\cdot 2^{m(N)}}+\frac{1}{10 m(N)}+
\eee
.$$
By \eqref{*e15}
\begin{equation}\label{*e25}
\bigg |A(f,x,N,m(N))-\frac{\sum_{k=\bbb_{m(N)-1}+1}^{N}f(T^{k}x)}
{N-\bbb_{m(N)-1}}\cdot \frac{\oN_{\bbb_{m(N)-1}}^{N}}{\oN_{0}^{N}}
\bigg |<
\end{equation}
$$\frac{2}{100\cdot 2^{m(N)}}+\frac{1}{10 m(N)}.$$

We also have
$$\bigg |\frac{\sum_{k=\bbb_{m(N)-1}+1}^{N}f(T^{k}x)}
{N-\bbb_{m(N)-1}}-\int f d\mu\bigg |=$$
$$\bigg |\frac{
N\frac{1}{N}\sum_{k=1}^{N}f(T^{k}x)-
\bbb_{m(N)-1}\frac{1}{\bbb_{m(N)-1}}
\sum_{k=1}^{\bbb_{m(N)-1}}f(T^{k}x)-
(N-\bbb_{m(N)-1})\int fd\mu}
{N-\bbb_{m(N)-1}}\bigg |\leq$$
$$\frac{|N\int fd\mu-\bbb_{m(N)-1}\int fd\mu-
(N-\bbb_{m(N)-1})\int f d\mu|+(N+\bbb_{m(N)-1})\eee}
{N-\bbb_{m(N)-1}}=$$
$$\frac{(N+\bbb_{m(N)-1})\eee}
{N-\bbb_{m(N)-1}}.$$
Using this in \eqref{*e25}
\begin{equation}\label{*e26}
\bigg |A(f,x,N,m(N))-
\frac{\oN_{\bbb_{m(N)-1}}^{N}}{\oN_{0}^{N}}
\int f d\mu\bigg |<
\end{equation}
$$\frac{2}{100\cdot 2^{m(N)}}+\frac{1}{10\cdot m(N)}+
\frac{\oN_{\bbb_{m(N)-1}}^{N}}{\oN_{0}^{N}}
\frac{N+\bbb_{m(N)-1}}{N-\bbb_{m(N)-1}}\eee.$$
Since $N\geq\bbb_{m(N)-1}$ we have $N+\bbb_{m(N)-1}\leq 2N$
and, obviously, $\oN_{\bbb_{m(N)-1}}^{N}/\oN_{0}^{N}\leq 1$.

To estimate $N-\bbb_{m(N)-1}$ we separate two
subcases.

CASE 1A. If $N-\bbb_{m(N)-1}>\sqrt\eee N$ then
$$\frac{N+\bbb_{m(N)-1}}{N-\bbb_{m(N)-1}}\eee<
\frac{2N}{\sqrt \eee N}\eee=2\sqrt \eee$$ and
from
\eqref{*e26}
it follows that
\begin{equation}\label{**e26}
\bigg |A(f,x,N,m(N))-\frac{\oN_{\bbb_{m(N)-1}}^{N}}{\oN_{0}^{N}}
\int fd\mu\bigg |<
\frac{1}{50\cdot 2^{m(N)}}+\frac{1}{10\cdot m(N)}+2\sqrt \eee.
\end{equation}

CASE 1B.
Suppose
 $N-\bbb_{m(N)-1}\leq \sqrt \eee N$. By \eqref{**f13} used
with $m_{0}=m(N)$ we have
\begin{equation}\label{***e26}
\frac{\oN_{\bbb_{m(N)-1}}^{N}}{\oN_{0}^{N}}
<
2\frac{N-\bbb_{m(N)-1}}{N}.
\end{equation}
Since $N-\bbb_{m(N)-1}<\sqrt\eee N$ we obtain
$$\frac{\oN_{\bbb_{m(N)-1}}^{N}}{\oN_{0}^{N}}< 2\sqrt\eee\text{
and }\frac{\oN_{\bbb_{m(N)-1}}^{N}}{\oN_{0}^{N}}\int
fd\mu< 2\sqrt\eee.$$
By its definition
$$A(f,x,N,m(N))=\frac{1}{\oN_{0}^{N}}\sum_{k=\oN_{m(N)-1}+1}^{\oN_{0}^{N}}
f(T^{n_{k}}x)\leq \frac{\oN_{\bbb_{m(N)-1}}^{N}}{\oN_{0}^{N}}
< 2\sqrt\eee,$$
and
\begin{equation}\label{*e27}
\bigg |A(f,x,N,m(N))-\frac{\oN_{\bbb_{m(N)-1}}^{N}}{\oN_{0}^{N}}
\int fd\mu\bigg |< 4\sqrt\eee.
\end{equation}
Therefore, in both cases (CASE 1A and CASE 1B) by \eqref{**e26},
or by \eqref{*e27} we have
\begin{equation}\label{*e28}
\bigg |A(f,x,N,m(N))-\frac{\oN_{\bbb_{m(N)-1}}^{N}}{\oN_{0}^{N}}
\int fd\mu\bigg |<\frac{1}{50\cdot 2^{m(N)}}+
\frac{1}{10m(N)}+4\sqrt\eee.
\end{equation}
Recalling \eqref{*d38}
we also have
\begin{equation}\label{**e28}
\sum_{m=1}^{m(N)-2}
A(f,x,N,m)<\frac{1}{100 m(N)},
\end{equation}
and we can suppose that our initial parameters were selected so
that
\begin{equation}\label{***e28}
\frac{\oN_{0}^{\bbb_{m(N)-2}}}{\oN_{0}^{N}}=
\frac{\oN_{m(N)-2}}{\oN_{0}^{N}}
\leq \frac{\oN_{m(N)-2}}{\oN_{m(N)-1}}
<\frac{1}{m(N)}.
\end{equation}
By using \eqref{**e24}, \eqref{*e28}, \eqref{**e28}, and
\eqref{***e28}
we conclude for CASE 1 that
$$|A(f,x,N)-\int f d\mu|\leq
|\sum_{m=1}^{m(N)-2}A(f,x,N,m)|+
|\frac{\oN_{0}^{\bbb_{m(N)-2}}}
{\oN_{0}^{N}}
\int fd\mu|+$$
$$|A(f,x,N,m(N)-1)-\frac{\oN_{\bbb_{m(N)-2}}^{\bbb_{m(N)-1}}}
{\oN_{0}^{N}}\int f d\mu|+
|A(f,x,N,m(N))-\frac{\oN_{\bbb_{m(N)-1}}^{N}}{\oN_{0}^{N}}\int fd\mu|<
$$
$$\frac{1}{100m(N)}+\frac{1}{m(N)}+\frac{2}{50\cdot 2^{m(N)}}+
\frac{2}{10m(N)}+4\sqrt\eee+\eee<5\sqrt\eee$$
if $N$ (and hence $m(N)$)
is sufficiently large (and $0<\eee<1$).
For CASE 2 from  \eqref{*e24} it also
follows that for large $N$'s we have
$|A(f,x,N)-\int f d\mu|<5\sqrt \eee.$

This implies that for any simple function
$0\leq f\leq 1$, and hence for an
arbitrary simple function
 the ergodic averages converge to the integral
of $f.$

\end{proof}

\end{document}